\documentclass[a4paper,oneside,12pt]{article}
\usepackage[left=2.5cm, right=2.5cm, top=2.5cm, bottom=3.5cm]{geometry}
\usepackage[english]{babel}
\usepackage[utf8]{inputenc}
\usepackage[T1]{fontenc}
\usepackage{authblk}
\usepackage{amsmath}
\usepackage{amsfonts}
\usepackage{amssymb}
\usepackage{amsthm}
\allowdisplaybreaks
\usepackage{enumitem}
	\setlist[itemize]{label={---}, itemsep=0pt, topsep=0pt, leftmargin=3em}
\usepackage{graphicx}
\usepackage{xcolor}
\usepackage{fancyhdr}
\usepackage{hyperref}

\setlist[enumerate]{label=(\roman*), itemsep=0pt, leftmargin=40pt}
\setlist[itemize]{label=---, itemsep=0pt, topsep=0pt, partopsep=0pt, leftmargin=40pt}
\setlist[description]{leftmargin=40pt}
 
\renewcommand{\le}{\leqslant}
\renewcommand{\ge}{\geqslant}
\newcommand{\pp}{\leqslant}
\newcommand{\pg}{\geqslant}

\newcommand{\EE}{\mathbb{E}}
\newcommand{\NN}{\mathbb{N}}
\newcommand{\PP}{\mathbb{P}}

\newcommand{\RR}{\mathbb{R}}

\newcommand{\ZZ}{\mathbb{Z}}

\newcommand{\cB}{\mathcal{B}}

\newcommand{\cD}{\mathcal{D}}

\newcommand{\cF}{\mathcal{F}}
\newcommand{\cG}{\mathcal{G}}
\newcommand{\cI}{\mathcal{I}}

\newcommand{\cS}{\mathcal{S}}
\newcommand{\cT}{\mathcal{T}}

\newcommand{\sE}{\mathsf{E}}
\newcommand{\sP}{\mathsf{P}}
\newcommand{\sQ}{\mathsf{Q}}

\newcommand{\petito}[1]{\mathrm{o}\left(#1\right)}
\newcommand{\grando}[1]{\mathrm{O}\left(#1\right)}

\DeclareMathOperator{\Div}{div}

\renewcommand{\limsup}{\varlimsup}
\renewcommand{\liminf}{\varliminf}
\newcommand{\dd}[1]{\,\mathrm{d}#1}
\newcommand{\Char}[1]{\mathbf{1}_{#1}}
\newcommand{\Ind}[1]{\mathbf{1}_{\left\lbrace #1 \right\rbrace}}

\newcommand{\tq}{\mathrel{/}}

\newcommand{\ens}[1]{\left\lbrace #1 \right\rbrace}
	
\newtheoremstyle{bfnote}%
{}{}%
{\itshape}{}%
{\sffamily\bfseries}{. ---\hspace{0.5em}}%
{ }%
{\thmname{#1}\thmnumber{ #2}\thmnote{ (#3)}}

\theoremstyle{bfnote}
\newtheorem{defin}{Définition}
\newtheorem{theorem}{Theorem}

\newtheorem{lemma}[theorem]{Lemma}

\begin{document}

\author[1]{Adrien \textsc{Perrel}}
\title{Limit~theorem for~subdiffusive random~walk in~Dirichlet random~environment in~dimension~$d \pg 3$}
\date{}
\affil[1]{Université Claude Bernard Lyon 1, \protect\\ Institut Camille Jordan, CNRS UMR 5208, \protect\\ Boulevard du 11 novembre 1918, 69622 Villeurbanne cedex, France \protect\\ (\href{mailto:perrel@math.univ-lyon1.fr}{\texttt{perrel@math.univ-lyon1.fr}}).}
\renewcommand\Affilfont{\small}

\maketitle

\begin{abstract}
We consider random walk in Dirichlet random environment in dimension ${d\ge 3}$, which corresponds to the case where the environment is constructed from i.i.d. transition probabilities at each vertex with a Dirichlet distribution with parameters $(\alpha_i)_{1 \pp i \pp 2d}$.
Dirichlet environments are weakly elliptic and the walk can be slowdowned by local traps whose strength are governed by a parameter $\kappa$. In this paper we prove a stable limit theorem when the walk is ballistic but subdiffusive, \emph{i.e.} when ${\kappa \in (1,2)}$. This completes the result of Poudevigne (\cite{poudevigne2019limit}) who proved a sub-ballistic stable limit theorem when $\kappa\in (0,1)$. Contrary to \cite{poudevigne2019limit}, we have to assume Sznitman's condition~$\mathbf{(T)}$ to prove the limit theorem since we work at the level of fluctuations and need a better control on renewal times.
\end{abstract}

\section*{Introduction}

This paper is concerned with a specific model of random walk in random environment on $\mathbb{Z}^d$, ${d\ge 3}$, with i.i.d. transition probabilities at each vertex (this model is referred to as RWRE in the sequel). As it is well-known, the non-reversibility of the environment is a major difficulty in the analysis of the asymptotic behavior of RWRE. While the one-dimensional case is very well understood, some fundamental questions are still open in higher dimension. The ballistic case, i.e. the case where an a priori ballistic condition is assumed, is by far the best understood. In this direction, Sznitman's conditions $\mathbf{(T)}$ or $(\mathbf{T})_{\gamma}$ (\cite{sznitman2001class, sznitman2002effective}), which ensure a good control on the displacement of the walk during two successive renewal times, play a central role. Under that condition and uniform ellipticity, annealed and quenched central limit theorems, large deviation thorems, have been obtained, (e.g. \cite{sznitman2001class, sznitman2002effective}, \cite{bolthausen2002static}, \cite{berger2008quenched},  \cite{rassoul2009almost}). It is still a major open question to prove that condition~$\mathbf{(T)}$ is equivalent to directional transience, even though important progresses have been made in that direction in proving equivalence of weaker and weaker forms of conditions~$\mathbf{(T)}$ (see \cite{campos2014ellipticity} and  \cite{bouchet2016ellipticity}). When the walk is not directional transient, the situation is far less clear, and the main results have been obtained in the perturbative regime (\cite{bolthausen2007multiscale}, \cite{baur2013long}, and \cite{sznitman2004diffusive} for a closely related model of diffusion in random environment), or for Dirichlet environments. 

In this work we investigate the slowdown effect of traps for random walk in random Dirichlet environment (RWDE) in the ballistic subdiffusive regime. The RWDE is a model of RWRE where the transition probabilities are independent at each site and are drawn according to Dirichlet random variables (see \cite{sabot2017random} for an overview). Due to the expression of the Dirichlet distribution, RWDE are not uniformly elliptic and the environment can create local traps where the walk spends an atypical large time. The Dirichlet environment is a good model to analyse this trapping effect for two main reasons. Firstly, it has been shown that Dirichlet environment have a key property of statistical invariance by time reversal for which one can infer a good control on exit probability of balls. This property has been used in several direction, e.g. to prove transience for $d\pg 3$, to prove the existence of an invariant measure viewed from the particle for $d\pg 3$, or to characterize directional transience (see \cite{tournier2009integrability}, \cite{sabot2011random}, \cite{sabot2013random} and \cite{tournier2015asymptotic}), and it is also instrumental in the present paper. Secondly, the traps are well understood and relatively simple to analyse (which is not the case for general weakly elliptic environment as shown in \cite{fribergh2016local}).

The local traps of Dirichlet environment are governed by an explicite parameter $\kappa$, which gives the tail exponent of the time spent in the strongest trap of the environment. In fact the strongest traps consist of simple edges on which the walk can backtrack a long time. In  \cite{poudevigne2019limit}, Poudevigne obtained stable limit theorems for the RWDE in the sub-ballistic case, i.e. for $\kappa\in (0,1)$, on $\mathbb{Z}^d$, $d\ge 3$. This paper completes this result and establishes stable fluctuations for RWDE in the ballistic sub-diffusive regime, i.e. $\kappa\in (1,2)$. But there is a major difference between our paper and \cite{poudevigne2019limit}. We need to assume condition~$\mathbf{(T)}$ in order to prove the result, while \cite{poudevigne2019limit} uses the existence of an absolutely continuous invariant measure viewed from the particle for an accelerated walk (already introduced in \cite{bouchet2013subballistic}). Correlatively, our approach is based on renewal times rather than on the accelerated walk as in \cite{poudevigne2019limit}. The reason why we need condition~$\mathbf{(T)}$ comes from the fact that we work at the level of fluctuations so we need to get a better control on renewal times of the accelerated walk than just the integrability given by the existence of an the invariant measure viewed from the particle. We borrow from \cite{poudevigne2019limit} several ideas as the partially forgotten walk or some estimates on time spent in traps. However, the fact that we need to control fluctuations, and not just the first order, leads to several extra difficulties.

Note that there is an important literature about slowdown for random walks in random environments for several families of models as Bouchaud trap model, biaised random walk on percolation clusters or in random conductances, random walks on Galton-Watson trees, 1D RWRE, we refer to Ben-Arous and Fribergh's lecture notes for an overview (\cite{benarous2016biased}). The difficulty in the model we consider comes from the highly non-reversible nature of the environment.

\section{Description of the model and known results.}

We restrict ourselves to the case of random walks in Dirichlet random environments to nearest neighbours in $\ZZ^d$, $d \pg 3$, but we will need to generalize this setting to more general directed graphs $\cG = (V, E)$ (see below section~\ref{prelim1}). In all cases, the graphs we consider are directed. Therefore, an edge $e \in E$ in $\cG$ is a couple $(x, y)$; we will denote $\underline{e} = x$ and $\overline{e} = y$. An unoriented edge $f \in \tilde{E}$ is by contrast a pair $\{x,y\}$; if $\{x, y\} \in \tilde{E}$; we will also denote $x \sim y$. Given $U$ a subset of vertices, we denote by $E_U$ and $\tilde{E}_U$ the sets of the respectively oriented and unoriented edges with both extremities in~$U$.

In this section, we only consider the euclidean lattice of $\ZZ^d$. We denote in the following $\{e_1, \ldots, e_d\}$ the canonical basis of $\ZZ^d$, as well as $e_{i+d}= -e_i$, so that $\{e_1, \ldots, e_{2d}\}$ is the set of unit vectors of $\ZZ^d$. The euclidean lattice is thus defined:
$$ V = \ZZ^d \text{\qquad and\qquad} E = \ens{(x,y)\in (\ZZ^d)^2 \tq y-x \in \{e_1, \ldots, e_{2d}\}}.$$

	\subsection{The random walk in random environment}
	\label{rwre_def}

\paragraph{Environments.} In this context, an admissible transition kernel on $\ZZ^d$ is an element of the set:
$$ \Omega = \ens{(\omega(x, y))_{(x,y) \in E} \in [0,1]^E \tq \forall x \in \ZZ^d, \sum_{y \tq (x,y)\in E} \omega(x,y) = 1}. $$
If we denote $\Upsilon = \ens{(p_i) \in [0,1]^{2d} \tq \sum_{i=1}^{2d} p_i = 1}$, we can also consider that $\Omega$ is the product set~$\Upsilon^{\ZZ^d}$. Endow $\Omega$ with the borelian $\sigma$-algebra~$\cB$. It is usual to call \emph{environment}, rather than transition kernel, the elements~$\omega$ of this set.

The model of RWRE consists in drawing at random an environment, or, in other words, in endowing $\Omega$ with a probability distribution~$\PP$.

\paragraph{Quenched law.}  Let us fix an environment~${\omega \in \Omega}$. We can now study the Markov chain~$X$ with transition kernel~$\omega$ issued from some initial probability distribution~$\mu$. This defines a probability distribution~$\sP^\omega_\mu$ over $V^\NN$ (endowed with the cylinder $\sigma$-algebra $\cF$) in the following way:
\begin{align*}
&\sP^\omega_\mu(X_0 = x) = \mu(x), \\
&\sP^\omega_\mu(X_{n+1} = y \,\vert\, X_n = x, X_{n-1}, \ldots, X_0) = \begin{cases}
\omega(x, y) &\text{if } (x,y) \in E \\
0 & \text{otherwise}
\end{cases}, \quad \forall n \in \NN,
\end{align*}
This probability distribution is called the quenched law of the RWRE in the environment~$\omega$, and we denote it $\sP^\omega_\mu$.

By definition, under $\sP^\omega$, $X$ is a Markov chain, and we can apply to~$X$ the Markov properties. For that, it is convenient to define the following stopping times: for every susbet $V \subset \ZZ^d$,
	\begin{align*}
	H_V &= \inf\ens{n \pg 0 \tq X_n \in V}; \\
	\bar H_V &= \inf\ens{n \pg 0 \tq X_n \not\in V}; \\
	H^+_V &= \inf\ens{n \pg 1 \tq X_{n-1} \not\in V, X_n \in V}.
	\end{align*}
In other words, $H_V$ is the first hitting time of~$V$; $\bar H_V$ is the first hitting time of the complementary of~$V$ (and therefore, if $X$ starts in~$V$, the first exit time from~$V$); $H_V^+$ is the first return time to $V$ if $X$ start in~$V$ (this supposes that $X$ has exited $V$ some time before returning in~$V$), or equal to $H_V$ otherwise.

\paragraph{Annealed law.} In order to take into account the two levels of randomness (choice of the environment, random walk in that environment), we define the annealed law $\sP_\mu$ on $\Omega \times V^\NN$ by:
	$$\forall (B, F) \in \cB \times \cF, \quad \sP_\mu(B \times F) = \int_B \sP_\mu^\omega(F) \, \dd{\PP(\omega)} = \EE\left[\sP_\mu^\omega(F) \Char{B}\right].$$

	\subsection{Dirichlet environment}
	\label{rwde_def}

In this work, we are interested in a peculiar probability distribution~$\PP$: the Dirichlet environment.

\paragraph{Dirichlet distribution.} The Dirichlet random environment derives its name from the Dirichlet distribution. Given $k \in \NN^*$ and $k$ positive parameters $(a_1, \ldots, a_k) \in (\RR_+^*)^k$, the Dirichlet distribution $\cD(a_1, \ldots, a_k)$ has density:
\begin{equation}
\label{eq:Dirichlet_distrib}
\frac{\Gamma\left(\sum_{i=1}^k a_i\right)}{\prod_{i=1}^{k} \Gamma(a_i)} \left(\prod_{i=1}^{k} x_i^{a_i-1}\right)
\end{equation}
	with respect to the Lebesgue measure over the simplex $\ens{(x_i)\in(\RR_+^*)^{k} \tq \sum_{i=1}^k x_i = 1}$.
	
	Alternatively, $\cD(a_1, \ldots, a_k)$ can be considered as the distribution of a random vector $\left(\frac{W_1}{\Sigma}, \ldots, \frac{W_k}{\Sigma}\right)$ where $(W_1, \ldots, W_k)$ are independent  random variables with respective law $\Gamma(a_i, \theta)$ and $\Sigma = \sum_{i=1}^k W_i$. It is therefore a natural generalization of the Beta distribution.

\paragraph{Dirichlet environment.} Fix $(\alpha_1, \ldots, \alpha_{2d})$ positive parameters. For $1 \pp i \pp 2d$, we may consider that every edge $(x,y)$ of the euclidean lattice such that $y-x = e_i$ is attributed a weight $\alpha_i$.

At each vertex $x \in \ZZ^d$, we draw independently at random a random vector $\omega_x = (\omega(x, x+e_i))_{1 \pp i \pp 2d}$ according to $\cD(\alpha_1, \ldots, \alpha_k)$.

	The joint distribution of the $\omega_x, x \in E$ is denoted $\PP$. When there is a risk of ambiguity, we may use the notation $\PP^{(\alpha)}$ and $\sP^{(\alpha)}$.
	
	As said in the introduction, Dirichlet environments have a specific property from which some asymptotic results of RWDE have been deduced. We enounce here the known results necessary to the comprehension of our main theorem. Other known results necessary to the proof are exposed in \ref{prelim1}, in particular the statistical invariance by time reversal.

\paragraph{Directional transience.} In Dirichlet random environments, much is known about the transience of the random walk. We first begin with directional transience:

\begin{theorem}[\cite{sabot2011random} and \cite{tournier2015asymptotic}]

Consider, in $\ZZ^d, d \pg 1$, the random walk in Dirichlet environment of parameter $(\alpha_1, \ldots, \alpha_{2d})$. Let us denote
$$
d_\alpha= \sum_{i=1}^{2d} \alpha_i e_i.
$$
Fix $\ell\in \RR^d\setminus\{0\}$.
\begin{enumerate}
	\item If $d_\alpha \cdot \ell =0$, then
\begin{align*}
	\sP^{(\alpha)}_0\text{-a.s.,}\qquad -\infty=\liminf_{n\to\infty} X_n\cdot \ell<\limsup_{n\to\infty} X_n\cdot \ell=+\infty. 
\end{align*}
	\item If $d_\alpha \cdot \ell >0$ (resp. $d_\alpha \cdot \ell < 0$)
\begin{align*}
	\sP^{(\alpha)}_0\text{-a.s.,}\qquad \lim_{n\to\infty}X_n\cdot \ell+\infty \qquad \text{(resp. } -\infty\text{).} 
\end{align*}
\end{enumerate}
Moreover, if  $d_\alpha\ne 0$, $(X_n)_n$ has an asymptotic direction given by $d_\alpha$: 
\begin{align*}
	P^{(\alpha)}_0\text{-a.s.,}\qquad \lim_{n\to\infty} \frac{X_n}{|X_n|}= \frac{d_\alpha}{|d_\alpha|}.
\end{align*}

\end{theorem}

\paragraph{Strongest traps in $\ZZ^d$.} In $\ZZ^d$ endowed with translation-invariant weights $(\alpha_1, \ldots, \alpha_{2d})$, if we want to investigate the asymptotic behaviour of the random walk in Dirichlet environment more precisely, we must find what are the strongest finite traps of the environment. Indeed, non uniform ellipticity of the environment can create traps of finite size, \emph{i.e.} finite subsets in which the random walk may spend an atypically large time.

We can measure how much traps delay the walk thanks to the order~$\kappa$ of the tail of the quenched Green functions $G^\omega(0, 0)$ of the walk. In the case of Dirichlet environment in dimension $d \pg 3$, Tournier proved in \cite{tournier2009integrability} that:

\begin{theorem}
\label{Green_integrabilite}
Assume $d\pg 3$. If $G_\omega(0,0)
$ denotes the Green function in environment $\omega$, then, for every $s \pp 0$:
	$$\sE^{(\alpha)}[G_\omega(0,0)^s]<\infty \qquad\Leftrightarrow \qquad s<\kappa.$$
where:
\begin{equation}
\label{eq:Kappa_Zd}
\kappa = \min_{1 \pp j \pp d} \kappa_j, \qquad \kappa_j = 2 \sum_{i=1}^{2d} \alpha_i - (\alpha_j + \alpha_{j+d}).
\end{equation}
\end{theorem}

\paragraph{Ballisticity regimes.} The parameter $\kappa$ also controls the ballisticity regime of the walk when $d\pg 3$ (see \cite{sabot2011random} and \cite{sabot2013random}).

\begin{theorem} Assume $d\pg 3$. Then, for all positive weights $(\alpha_1, \ldots, \alpha_{2d})$,
$$ | X_n | \xrightarrow[n\to\infty]{\sP_0^{(\alpha)}\text{-a.s.}} +\infty. $$

Moreover, set $d_\alpha = \sum_{i=1}^{2d} \alpha_i e_i$; then:

(i) If $d_\alpha \neq 0$ and $\kappa >1$, there exists $v\in \RR^d\setminus\{0\}$ such that
$$
\lim_{n\to \infty} {X_n \over n} = v.
$$

(ii) If $d_\alpha = 0$ or $\kappa\pp 1$, then
$$
\lim_{n\to \infty} {X_n \over n} = 0.
$$
\end{theorem}

	\subsection{The (T) condition}

We aim in this paper at the generalization of the results of Poudevigne \cite{poudevigne2019limit} in the subdiffusive case, i.e. when $\kappa \in (1, 2)$. We enounce rigorously Poudevigne's result in Theorem \ref{poudevigne}, at the end of this section. Poudevigne finds the limit law of the process $(n^{-\kappa} X_{\lfloor nt \rfloor})_{t \pg 0}$ when $n \to \infty$, that is to say informally, the limit law of the fluctuations of $(X_n)_{n \in \NN}$ at order $n^{\kappa}$, in the subballistic case, i.e. when $\kappa < 1$.

However, a new difficulty appears when $\kappa > 1$.

When $\kappa < 1$, the fluctuations in time are of order $n^{1/\kappa}$ at time~$n$, and thus are greater than the typical distance done at time~$n$, which is of order $n^\kappa$ (see \cite{bouchet2013subballistic}).

In our case $\kappa > 1$, we have to make a distinction between the typical distance, of order $n$ at time~$n$, and the fluctuations in space and in time, which are expected to be of order $n^{1/\kappa}$, lesser than $n$. The main problem in controlling these fluctuations is that we do not know if the time period between two renewal times (see definition below) has finite $\gamma$-moment for some $\gamma > 1$.

So it becomes necessary to make further assumptions. The most common one in the study of RWRE is condition~$\mathbf{(T)}$, introduced in \cite{sznitman1999law}.
	
\paragraph{Renewal structure.}
Let $\hat{u}$ be a unit vector of $\RR^d$. We define $(T_i)_{i\in\NN^*}$, the renewal times in the direction $\hat{u}$, as follows.

Denote by $(\Theta_n)_{n \pg 1}$ the canonical time shifts on
$( \mathbb{Z}^d )^\mathbb{N}$.
For $l\geq 0$, we define the time
\[
\Delta_l = \min\{ n \pg 0 \tq \langle X_n, \hat u \rangle \geq l \}.
\]

Fix $a > 2 \sqrt{d}$ (this condition is only necessary to ensure that $\ens{x \in \ZZ^d \tq 0 \pp \langle x, \hat u \rangle \pp a}$ is a connected set). We define, for initialization:
\begin{align*}
D &= \min\{ n \geq0 \tq \langle X_n , \hat u \rangle < X_0 \}, \\
M_0 &= \langle X_0, \hat u \rangle \\
\tau_0 &= 0;
\end{align*}
and then, recursively for $k \geq 0$,
\begin{align*}
\tau_{k+1} &= \Delta_{M_k+a}, \\
\tau'_{k+1} &= D \circ \Theta_{\tau_{k+1}} + \tau_{k+1}, \\
M_{k+1} &= \max\ens{ \langle X_n, \hat u \rangle, 0 \leq n \leq \tau'_{k+1} }.
\end{align*}

The first regeneration time is then defined as
\[
T_1 = \tau_K \qquad\text{ where } K = \min\{ k \geq 1 \tq \tau_k < \infty,
\tau'_k = \infty\}.
\]
We can now define recursively in $n$ the $(n + 1)$-th regeneration time $ T_{n+1} $ as:
$$ T_{n+1} = T_1 (X_\cdot) + T_n ( X_{T
_1+\cdot} - X _{T_1} ) .$$
We may in the following write $T_k(\hat u)$ to emphasize the dependence on the direction~$\hat u$.

Prop. 1.2 of \cite{sznitman1999law} shows that, if $X$ is transient in direction $\hat u$, then $\sP_0(D = 0) > 0$, and that $T_1$ is finite $\sP_0$-almost surely.

The renewal times are used to create independence:

\begin{lemma}[Prop. 1.3 of \cite{sznitman1999law}]
For all $k\in\NN^*$, let $\mathcal{G}_k$ be the $\sigma$-field defined by:
\[
\mathcal{G}_k:=\sigma(T_1,\dots,T_k,(X_n)_{0\leq n \leq T_k},(\omega(x,\cdot))_{\langle x, \hat u \rangle < \langle X_{T_k}, \hat u \rangle}),
\]
and set $\sQ_0 = \sP_0(\cdot \vert T_1 < \infty).$
We have, for all $k\geq 1$:
\begin{multline*}
\sQ_0\left((X_{T_k+n})_{n\geq 0} \in \cdot, (\omega(X_{T_k}+x,\cdot))_{\langle x, \hat u \rangle \geq 0}\in \cdot |\mathcal{G}_k\right)\\
=\sP_0\left((X_n)_{n\geq 0} \in \cdot, (\omega(x,\cdot))_{\langle x, \hat u \rangle \geq 0}\in \cdot |D=0\right).
\end{multline*}
\end{lemma}
This means that the trajectories and the transition probabilities inside slabs between two consecutive renewal times (after the first one) are i.i.d random variables.

\paragraph{The (T) condition.} In this paper, we constantly assume condition~$\mathbf{(T)}$, which we recall below.

Given a unitary vector $\hat{u}$, we say that the condition~$\mathbf{(T)}$ in direction $\hat{u}$ is satisfied if:
	\begin{enumerate}
	\item $X$ is transient in the direction~$\hat u$, i.e.:

\begin{equation}
\label{eq:Ti} \tag{\textbf{T}i}
\langle X_n, \hat u \rangle \xrightarrow[n\to\infty]{\sP_0-\text{a.s.}} +\infty;
\end{equation}
	\item There exists $c > 0$ such that:
\begin{equation}
\label{eq:Tii} \tag{\textbf{T}ii}
\sE_0\left[\sup_{1 \le i \le T_1} e^{c \|X_i\|} \right] < \infty.
\end{equation}
	\end{enumerate}

In the case of Dirichlet random environments, Tournier have shown in \cite{tournier2009integrability} that $\mathbf{(T)}$ is satisfied as soon as
$$ \sum_{i=1}^d \vert \alpha_i - \alpha_{i+d} \vert > 1. $$
It is thus possible to find parameters $\alpha$ such that $\kappa \in (1, 2)$ and $\mathbf{(T)}$ is satisfied.

In fact, it is conjectured (but this conjecture is far from being proved) that $\mathbf{(T)}$ is equivalent to the directional transience of the walk; in that case, assuming condition~$\mathbf{(T)}$ would not be a restriction of generality in our paper.

\paragraph{Subballistic regime.} When $\kappa < 1$, Poudevigne has established in \cite{poudevigne2019limit} the law of the fluctuations of the random walk. Remind that the L\'evy process $\mathcal{S}^{\kappa}$ with completely asymmetric $\kappa$-stable increments is characterized by the following generator function:
\[
\forall \lambda\in\RR,\forall s\in \RR^+, \EE\left(\exp\left(i\lambda \mathcal{S}^{\kappa}_s\right)\right)=\exp\left(-s|\lambda |^{\kappa}\left(1-i\text{sgn}(\lambda)\tan\left(\frac{\pi\kappa}{2}\right)\right)\right)
\]
When $\kappa \in (0,1)$, since this process is non-decreasing and c\`adl\`ag, its c\`adl\`ag inverse $\tilde{\mathcal{S}}^{\kappa}$ can be defined as follows:
\[
\tilde{\mathcal{S}}^{\kappa}_t:= \inf\{s, \mathcal{S}^{\kappa}_s\geq t\}.
\]

\begin{theorem}[Theorem~1 of \cite{poudevigne2019limit}]
\label{poudevigne}
Set $d\geq 3$ and $(\alpha_1, \ldots, \alpha_{2d})\in (\RR_+^*)^{2d}$. Let $X^n(t)$ be defined by:
\[
X^n(t)=n^{-\kappa} X_{\lfloor nt \rfloor}.
\]
If $\kappa<1$ and $\sum_{i=1}^{2d} \alpha_i e_i \neq 0$, there exists positive constants $c_1, c_2,c_3$ such that in the $J_1$-topology:
\[
\left(t\longmapsto n^{-\frac{1}{\kappa}}T_{\lfloor nt \rfloor}\right)\xrightarrow[n\to\infty]{\sP_0^{(\alpha)}\text{-(d)}} c_1 \mathcal{S}^{\kappa};
\]
in the $M_1$-topology:
\[
\left(t\longmapsto n^{-\frac{1}{\kappa}}\inf\{t\pg 0, X^n(t).e_1\pg nx\}\right) \xrightarrow[n\to\infty]{\sP_0^{(\alpha)}\text{-(d)}} c_2 \mathcal{S}^{\kappa};
\]
and in $J_1$-topology:
\[
X^n \xrightarrow[\to\infty]{\sP_0^{(\alpha)}\text{-(d)}} c_3 \tilde{\mathcal{S}}^{\kappa}d_{\alpha}.
\]
\end{theorem}
We do not discuss the different topologies of the space of càdlàg functions: for details, see \cite{whitt2002limits}.

We now are ready to state our main theorem concerning the subdiffusive regime, and compare it with the subballistic regime.

\section{Main results}

	\subsection{Limit theorem}
	
The main result of this paper is the following.

\begin{theorem}
\label{main}
Assume that $d \pg 3$, $\kappa \in (1, 2)$ and $d_\alpha \neq 0$ and condition $(T)$ in some direction $\hat{u}$. Let $X$ be a random walk in a Dirichlet environment, and define
$$ S_t = \sup\{\langle X_k, \hat{u} \rangle, k \pp t\}. $$

Then there exists constants $c_1, c_2, c_3 > 0$, $\tau = \sE_0^{(\alpha)}[T_2-T_1]$ and $v = \frac{\sE_0^{(\alpha)}[X_{T_2} - X_{T_1}]}{\sE_0^{(\alpha)}[T_2 - T_1]}$ such that, in the $J_1$-topology:
\[
\left(t \longmapsto n^{-\frac{1}{\kappa}}\left(T_{\lfloor nt \rfloor} - n\tau t \right)\right)\xrightarrow[n\to\infty]{\sP_0^{(\alpha)}\text{-(d)}} c_1 \mathcal{S}^{\kappa}
\]
in the $M_1$-topology:
\[
\left(t\longmapsto n^{-\frac{1}{\kappa}} \left(S_{nt} - n \langle v, \hat u \rangle t\right) \right) \xrightarrow[n\to\infty]{\sP_0^{(\alpha)}\text{-(d)}} c_2 \mathcal{S}^{\kappa};
\]
and in the $M_1$-topology:
\[
\left(t\longmapsto n^{-\frac{1}{\kappa}} \left( X_{\lfloor nt \rfloor} - nvt \right) \right) \xrightarrow[n\to\infty]{\sP_0^{(\alpha)}\text{-(d)}} -c_3 {\mathcal{S}}^{\kappa} d_{\alpha},
\]
where $\cS^\kappa$ is the $\kappa$-stable Lévy process.
\end{theorem}

We have already discussed the necessity of the $\mathbf{(T)}$-condition above. This condition implies in a way that fluctuations are essentially time-related (and not space-related): and indeed, there are no fluctuations orthogonal to the asymptotic speed. Therefore, the sources of the fluctuations are likely to be the traps, which delay the walk without having a big impact on its trajectory. We thus define precisely this notion of trap, before giving a sketch of our proof.

	\subsection{Notion of trap}
	\label{traps}
	
	It is known that, because of non-uniform ellipticity of the environment, traps of finite size are created, \emph{i.e.} (in the general meaning of the term) finite subsets in which the random walk spends an atypically large time. Nonetheless, the precise definition of trap varies according to the authors.
	
	To define specifically traps in this article, we will follow~\cite{poudevigne2019limit}. In particular, we consider that traps can only be undirected edges, and not larger subsets of~$\ZZ^d$. That is enabled by Theorem~\ref{Green_integrabilite}: for further explanation, see the discussion following Theorem~\ref{Green_integr} below. 

\paragraph{Strength of an undirected edge.} Given an undirected edge $f = \{x, y\}$, we define its back-and-forth probability
\begin{equation}
\label{eq:pf}
p_f = \omega(x,y) \omega(y,x)
\end{equation}
and its strength
\begin{equation}
\label{eq:sf}
s_f = \frac{1}{2 - \omega(x,y) - \omega(y,x)}.
\end{equation}
We also denote, for every $x \in \ZZ^d$,
$$ s_x = \max_{f \ni x} s_f. $$

It readily follows from these definitions that, for every undirected edge~$f$:
\begin{equation}
\label{eq:pf-sf-general}
(1-p_f) s_f \pp 1.
\end{equation}

\paragraph{Traps.} A \emph{trap} is any undirected edge $f = \{x,y\}$ such that $\omega(x,y)+\omega(y,x)> \frac{3}{2}$. Note that, if $f = \{x, y\}$ is a trap, $\omega(x,y) > \frac{1}{2}$ and $\omega(y,x)>\frac{1}{2}$.

As a first consequence, for every point $x \in \ZZ^d$, there is at most one neighbour $y$ such that $\{x,y\}$ is a trap.

As a second consequence, when $f = \{x, y\}$ is a trap, then:
\begin{align*}
1 - p_f &= (1-\omega(x,y)) + (1-\omega(y,x)) - (1-\omega(x,y)) (1-\omega(y,x)) \\
&\pg \frac{1}{2} (1-\omega(x,y)) + (1-\omega(y,x)) - \frac{1}{2} (1-\omega(y,x)) \\
&\pg \frac{1}{2s_f},
\end{align*}
so that, when $f$ is a trap:
\begin{equation}
\label{eq:pf-sf-trap}
(1-p_f) s_f \pg \frac{1}{2}.
\end{equation}

Let $\mathcal{T}^{\omega}$ (or~$\mathcal{T}$ when there is no ambiguity) be the set of traps~$\{x,y\}\in \tilde{E}$ for the environment~$\omega$, and $\tilde{\mathcal{T}}^{\omega}$ (or~$\tilde{\mathcal{T}}$) the set of the extremities of these traps. We will add an index $U$ when refering to the traps~$\mathcal{T}^{\omega}_U$ contained in~$U$ and to their extremities~$\tilde{\mathcal{T}}^{\omega}_U$, and an index~1 when refering to the traps~$\mathcal{T}^{\omega}_1$ encountered by the walk before~$T_1$ and to their extremities~$\tilde{\mathcal{T}}^{\omega}_1$.

\paragraph{Time spent in traps.} In the following, we will make use of the following durations:
\begin{itemize}
\item $T_1^\circ$ the time spent before $T_1$ (excluded) outside of traps:
$$ T_1^\circ = \sum_{i=1}^{T_1} \Ind{X_{i-1} \in V\setminus\tilde{\cT}}; $$
\item when $f$ is an undirected edge, $N_f$ is the number of visits of the walk to this edge (that is to say, the number of times it enters this edge, or equivalently the number of time it exits this edge) before time $T_1$:
$$ N_f = \sum_{i= 0}^{T_1-1} \Ind{X_i \in f, X_{i+1} \not\in f}; $$
\item when $f$ is an undirected edge and $1 \pp j \pp N_f$, $\ell_f^j$ is the number of steps the walk takes in $f$ during its $j$-th visit: for example
$$ \ell_f^1 = \sum_{i \pg 0} \Ind{\forall 0 \pp k \pp i, X_{H_f + k + 1} \in f}; $$
\item $\cT_1^\omega$ is the set of traps visited by the walk before time $T_1$ (excluded) and $T_1^\bullet$ is the time spent before $T_1$ (excluded) in traps:
$$ T_1^\bullet = \sum_{i=1}^{T_1} \Ind{X_{i-1} \in \tilde{\cT}} = \sum_{f \in \cT_1} \left(N_f + \sum_{j=1}^{N_f} \ell_f^j\right). $$
\end{itemize}
	
	\subsection{Sketch of proof}
	\label{sketch}
	
As explained earlier, in order to prove the first clause of the theorem, we follow a strategy similar to \cite{kesten1975limit} we estimate the tail of $T_1$, and then apply the stable law theorem (and the renewal structure) to obtain the fluctuations of $T_n$. We have already announced that the main contribution in this tail is the time spent in traps. We must so first show that the fluctuations of the time spent outside traps is negligible.
	
As for the traps, \cite{poudevigne2019limit} has showed that the order $n^{1/\kappa}$ of the fluctuations is due to the time spent in traps of optimal direction: informally, knowing the trajectory of the walk out of traps, the tail of the time spent in a trap colinear to the direction $e_j$ decreases as $C_\mathfrak{c} x^{-\kappa_j}$, with a multiplicative constant $C_\mathfrak{c}$ that only depends on the configuration $\mathfrak{c}$ of the trap, which is a random variable which essentially retains the number of visits to the trap and the vertices by which the walk enters and exits the trap. To draw all the conclusions of this estimate---that is to say, informally, to integrate this tail bound,--- it is necessary to control $C_\mathfrak{c}$. That leads us to separate two cases: when the number of visits of the traps is high, and when it is moderate. We so introduce a threshold $m: \epsilon \in \RR_+ \longrightarrow m(\epsilon) \in \NN$, which is to be sent to $\infty$ by letting $\epsilon \to 0$ at the end of the proof, that separates the traps visited a few times and those visited often.

So we are led to this decomposition of $T_1$: for every $\epsilon > 0$,
$$ T_1 = T_1^\circ + \sum_{f \in \cT_1} N_f + \sum_{f \in \cT_1}\sum_{j=1}^{N_f} \ell_f^j \Ind{N_f \pg m(\epsilon)} + \sum_{f \in \cT_1}\sum_{j=1}^{N_f} \ell_f^j \Ind{N_f \pp m(\epsilon)}. $$

Every term of that sum ($T_1^\circ$, $N_f$, $\ell_f^j$) can be controlled; however, the number of terms of the sum, \emph{i.e.} the cardinality of $\cT_1$, is hard to control because it is closely related to the number of points the walk visits before $T_1$, or, in other terms, it is closely related to the spatial behaviour of the walk. Here appears the necessity to assume $\mathbf{(T)}$: it enables us to assume that, if $T_1 \pg x$, with very high probability (let us say of order at least $1 - x^{-\kappa-\epsilon}$), the trajectory walk is contained in a box $U$ of side $\ell(x)$ of order $\log(x)$.

Finally, we follow the following strategy:
	
\begin{enumerate}
\item We show that, with high probability, we can assume that, before $T_1$, the random walk stays confined in a box $U$ whose edge has length of order $\ln(x)$.
\item We estimate the moments of the time spent out of traps, $T_1^\circ$, or equivalently, of the sum of the number of visits to the vertices which do not belong to a trap. The condition of not belonging to a trap can be controlled by Markov's inequality, while the estimation of the moments of the number of visits to a vertex is easily rephrased, by common arguments about Green functions, in the problem of the estimation of return probabilities. In the Dirichlet case, these return probabilities can be estimated by using the property of statistical invariance by time reversal of this environment (see \cite{sabot2017random} for introductory examples).
\item The number of visits of traps can be approximated by the number of visits to a vertex in a contracted graph. The impact of contraction on Dirichlet environment was first studied in \cite{tournier2009integrability}. We use an analogous strategy to compare the moments of the number of visits to traps with the Green functions in the contracted graph.
\item The total time spent in a trap visited more than $m(\epsilon)$ times is roughly a sum of at most $m(\epsilon)$ geometric random variables whose (random) parameter is determined by the strength of the trap. In that sum, we can first neglect traps whose strength is greater than a threshold $h(x)$ of order $\epsilon x$, because it is unlikely for the walk to visit them. We can then focus on traps of strength less than $h(x)$. For these ones, the approximation by a sum of geometric variables with bounded parameter is sufficient to show that the tail of the sum of the number of visits to them is negligible.
\item As for the traps visited less than $m(\epsilon)$ times, in a similar fashion as in the previous point, we can first neglect those whose strength is less than $h(x)$, and then we focus on those with strength bigger than $h(x)$. We show that, with high probability, there is only one single such trap (that is: visited less than $m(\epsilon)$ times and with strength bigger than $h(x)$). For this one, we can use the tail estimate of the time spent in a single fixed trap established in \cite{poudevigne2019limit}: this tail is of order $C(\epsilon) x^{-\kappa}$.
\item We let $\epsilon$ go to~0 (\textit{i.e.} $m(\epsilon)$ go to~$\infty$): the predominant term is the one of order $C(\epsilon) x^{-\kappa}$. Because $C(\epsilon)$ grows as $\epsilon$ goes to~0, it tends to some (possibly infinite) constant $C_0$; but $C_0 - C(\epsilon)$ being linked to all the negligible partial tails we have enumerated above, it cannot be infinite. This establishes that $T_1$ has tail of order $C_0 x^{-\kappa}$.
\end{enumerate}

The tail estimate implies, by a stable law theorem, a limit theorem for $T_n$. 

Finally, to deduce from the limit theorem for $T_n$ a limit theorem for $X_n$, we use common inversion techniques, exposed for example in \cite{whitt2002limits}.
	
\section{Limit theorem for renewal times}

	We fix a large $x > 0$ and a parameter $h(x)$ to be ajusted later on. Because of the renewal structure of the random walk, the main goal of this proof is the estimation of $\PP(T_1 > x | B)$, where $B$ is the event that the walk never steps in the halfspace $\ens{(x_1, \ldots, x_d) \in \ZZ^d \tq x_1 < 0}$.
	
	As explained in the sketch of proof, we decompose the duration $T_1$ in five parts:
	\begin{align*}
	T_1 = T_1^\circ &+ \sum_{f \in \cT_1} N_f + \sum_{f \in \cT_1}\sum_{j=1}^{N_f} \ell_f^j \Ind{s_f \pp h(x)} \\*
	&+ \sum_{f \in \cT_1}\sum_{j=1}^{N_f} \ell_f^j \Ind{s_f > h(x), N_f \pp m(\epsilon)} + \sum_{f \in \cT_1}\sum_{j=1}^{N_f} \ell_f^j \Ind{s_f > h(x), N_f > m(\epsilon)}.
\end{align*}
and estimate the tail of each term of this sum separately.

Moreover, we also explained that we will first work under the further assumption that $T_1 < \bar H_U$, where $U$ is a fixed box $U = [0, \ell(x)] \times [-\ell(x), \ell(x)]^{d-1}$. By condition $\mathbf{(T)}$, there exists a constant $c_T > 0$ such that:
$$ \PP(T_1 > \bar H_U) \pp e^{-c_T \ell(x)}. $$
Since we expect $T_1$ to have a tail of order $x^{-\kappa}$, we choose:
\begin{equation}
\label{eq:l(x)}
\ell(x) = \frac{1+\kappa}{c_T} \log(x)
\end{equation}
so that, for $x \to \infty$:
$$ \PP(T_1 > \bar H_U) = \textrm{o}(x^{-\kappa}). $$

	\subsection{Preliminary 1: Dirichlet environment on general graphs}
	\label{prelim1}
	
We will be compelled later to work on graphs different from $\ZZ^d$, for example finite subgraphs of $\ZZ^d$ or contracted graphs deduced from subgraphs. We generalize in this section the Dirichlet setting to these new graphs.

\paragraph{Dirichlet environment.} Let be given an oriented graph $\cG = (V, E)$, whose edges~$e$ are endowed with a weight $\alpha(e) > 0$. We draw independently at random, for every vertex~$x\in V$, the random variable $\omega_x = (\omega_e)_{\underline{e}=x}$ according to the Dirichlet distribution $\cD\left(\ens{\alpha(e), \underline{e}=x}\right)$ (see equation~\ref{eq:Dirichlet_distrib} for the definition of a Dirichlet random variable). This endows $\cG$ with an environment $(\omega_x)_{x\in V}$, whose distribution is called Dirichlet environment with parameters $\ens{\alpha(e), e \in E}$ and is denoted~$\PP^{(\alpha)}$.
	
\paragraph{Change of measure.}	In particular, on a finite graph $\cG$ endowed with weights $\alpha(e), e \in E$, we have, for every $(\xi(e))\in \RR^E$ such that $\alpha(e)+\xi(e)>0$:

\begin{equation}
\label{eq:changement_mesure}
\frac{\dd \PP^{(\alpha)}}{\dd \PP^{(\alpha+\theta)}} = \frac{Z_\alpha}{Z_{\alpha+\xi}} \prod_{e\in E} \omega(e)^{-\xi(e)}.
\end{equation}
where, for every function $(\theta(e))\in \RR^E$, we denote
$$ Z_\theta = \frac{\prod_{e\in E} \Gamma(\theta(e))}{\prod_{x\in V} \Gamma(\theta(x))}. $$
	
	If $\omega$ is drawn at random following $\PP^{(\alpha)}$, its joint moments are as follows: for every $(\xi(e))\in \RR^E$ such that $\alpha(e)+\xi(e)>0$ for all $e \in E$, 
\begin{align*}
	\EE^{(\alpha)}\left[ \prod_{e\in E} \omega(e)^{\xi(e)}\right]
= \left(\prod_{e\in E} \frac{\Gamma(\alpha(e)+\xi(e))}{\Gamma(\alpha(e))}\right)
\left(\prod_{x\in V} \frac{\Gamma(\alpha(x))}{ \Gamma(\alpha(x)+\xi(x))}\right) = \frac{Z_{\alpha+\xi}}{Z_\alpha},
\end{align*}
where, for every function $(\theta(e))\in \RR^E$, we have denoted
$$ \theta(x)=\sum_{e \in E \tq \underline e =x} \theta(e), \qquad Z_\theta = \frac{\prod_{e\in E} \Gamma(\theta(e))}{\prod_{x\in V} \Gamma(\theta(x))}. $$
When $\xi_e+\alpha_e \le 0$ for some edge $e$, the expectation is infinite.
	
	\paragraph{Divergence and time reversal.} For every vertex~$x \in V$, we call \emph{divergence} of~$\alpha$ at~$x$ the difference between the weights going in and going out from this vertex:
	$$ \Div(\alpha)(x) = \sum_{(x,y)\in E} \alpha(x,y) - \sum_{(y,x)\in E} \alpha(y,x). $$
	We will always deal with graph with null divergence in the following.

Consider a finite oriented graph $\cG = (V,E)$ endowed with positive weights $(\alpha(e))_{e\in E}$; and define its reversed graph $\check{\cG} = (V, \check{E})$ with same vertices but reverted edges: for every $(x, y)\in V^2$, $e = (x,y) \in E$ if and only if $\check{e} = (y,x) \in \check{E}$, and in this case $\check{e}$ is endowed with a weight $\check{\alpha}(\check{e}) = \alpha(e)$. Now, assume that $\cG$ is endowed with an environment~$\omega \in \Omega$. Since $\cG$ is finite, there exists an invariant measure $\pi^\omega$ over $\cG$ with respect to $\omega$. We then define over $\check{\cG}$ a new environment $\check\omega$ by:
$$\forall (x,y) \in E, \quad \check\omega(y,x) = \frac{\pi^\omega(x)}{\pi^\omega(y)} \omega(x,y). $$

\begin{lemma}[Time reversal property]
Assume that $\cG$ is finite and that $\Div(\alpha) \equiv 0$. Then:
$$ \omega \sim \sP^{(\alpha)} \Rightarrow \check{\omega} \sim \sP^{(\check{\alpha})}. $$
\end{lemma}

The first important consequence of the time reversal property is the following:

\begin{lemma}
Assume that $\cG = (V, E)$ is finite and that $\Div(\alpha) \equiv 0$. Then, for every $(x, y)\in V^2$,
$$\sP_x^\omega(H_y > H_x^+) \pg \sP_x^{\check{\omega}}(X_1 = y);$$
and under $\sP^{(\alpha)}$, this last probability follows a Beta distribution with parameters $(\alpha(y,x), \alpha_x)$.
\end{lemma}

This result is the main ingredient in the study of transience in the Dirichlet environment over $\ZZ^d$ (we refer to \cite{sabot2004ballistic}, with a summary in \cite{sabot2017random}). If we apply it to a finite graph of diameter $n$ properly constructed (take a ball around $0$ in $\ZZ^d$, add a vertex $\partial$ out of this ball, link every vertex on the border of this ball to $\partial$ and $\partial$ to 0) and endowed with the weights naturally induced by those of $\ZZ^d$ and completed so as to have null divergence, we can estimate the probability starting from~$0$ to reach a point at distance~$n$ before returning to~$0$, and take the limit $n \to \infty$ to obtain information about transience.
	
\paragraph{The kappa parameter.} We now aim at generalizing the $\kappa$ parameter introduced for $\ZZ^d$ in equation~\ref{eq:Kappa_Zd}.

For every finite subset of vertices $S$ of $\cG$, we define its $\kappa$-parameter:
\begin{equation}
\label{eq:Kappa_S}
\kappa(S) = \min \ens{ \sum _{e \in \partial_+(K)}  \alpha _e; K \text{ connected set of vertices},\, \{0\} \subsetneq K \text{ and } \partial S \cap K \neq \emptyset }
\end{equation}
where we denote
$$\partial S = \ens{x \in S \tq \exists y \not\in S, x \sim y}, \qquad \partial_+ S = \ens{e \in E \tq \underline{e} \in S, \overline{e} \not\in S}. $$
This $\kappa(S)$ exactly measures how much $S$ delays the walk. Indeed, if we denote $G_S^\omega(x, x)^s$ the quenched Green function of the walk in~$S$ starting from~$x \in S$, we have the following result, proved by Tournier in~\cite{tournier2009integrability}:

\begin{theorem}
\label{Green_integr}
Let $\cG = (V \cup \{\partial\}, E)$ be a finite oriented graph such that every vertex of $V$ is connected to $\partial$ in $\cG$ and such that for every $e\in E$, $\underline{e} \neq \partial$. Endow $\cG$ with a Dirichlet environment with positive weights $(\alpha(e))_{e\in E}$. Let $x \in V$.

Let $S$ be a finite subset of~$V$. Then, for every $s \pp 0$ and $x \in S$,
$$\sE^{(\alpha)} \left[G_S^\omega(x, x)^s\right] < \infty \qquad \Leftrightarrow \qquad s < \kappa(S).$$
\end{theorem}

With definition~\ref{eq:Kappa_S}, in $\ZZ^d$, the $\kappa$-parameters introduced in equation~\ref{eq:Kappa_Zd} could also be reformulated as follows:
\begin{align*}
\kappa_j &= \kappa(\ens{0, e_j}), \\
\kappa &= \min \ens{\kappa(S), S \text{ finite connected set of vertices of } \ZZ^d} 
\end{align*}
In other words, in $\ZZ^d$, $\kappa$ is the value of the minimal $\kappa$-parameter of finite subsets, and this value, being one of the $\kappa_j$'s, is reached by some pairs of contiguous vertices---or, in other words, by some undirected edges. Comparing theorems~\ref{Green_integrabilite} and~\ref{Green_integr}, they suggest that, in~$\ZZ^d$, there is not any infinite subset of vertices than can have a bigger impact on the walk than undirected edges, since the global quenched Green function is as integrable as the quenched Green function restricted to the strongest finite trap. Undirected edges can therefore be considered as the strongest finite traps (in the general meaning) of~$\ZZ^d$, and this explains why, in section~\ref{traps}, we considered that traps (in the specific meaning of this article) could only be undirected edges. The strength~$s_f$ of an undirected edge (defined in equation~\ref{eq:sf}) can be considered as a ``quenched'' measure of its strength, while $\kappa(f)$ is an ``annealed" one. The link between these two measures of strength is that $s_f$ has tails of order $x^{-\kappa(f)}$ (see below lemma~\ref{quasi-independence}.

However, for example when considering the time spent by the walk outside of traps (in the specific meaning of this article), we will have to measure the strength of larger subsets of vertices of~$\ZZ^d$ than undirected edges. We thus define:
\begin{align*}
\kappa' = \min \ens{\kappa(S), S \text{ connected set of 3 vertices of } \ZZ^d}.
\end{align*}
In particular, as stated in the discussion of Theorem~\ref{Green_integr}, $\kappa' > \kappa$.

	\subsection{Preliminary 2: partially forgotten walk}
	\label{prelim2}

	We determine as a first step an equivalent of the tail of the time spent in a given trap. It gives the principal order of the fluctuations.
	
\paragraph{Trap-equivalent environments.} We use the definitions of \cite{poudevigne2019limit}.

We say that two environments $\omega_1$ and $\omega_2$ are trap-equivalent if:
\begin{itemize}
\item they have the same traps: $\mathcal{T}^{\omega_1}=\mathcal{T}^{\omega_2}$;
\item at each vertex not in a trap, the transition probabilities are the same for both environment:
\[
\forall x \not\in\tilde{\mathcal{T}}^{\omega_1},\ \forall y\sim x,\ \omega_1(x,y)=\omega_2(x,y),
\]
\item at each vertex $x$ in a trap $\{x,y\}$, the transition probabilities conditioned on not crossing the trap are the same:
\[
\forall (x,y)\in E,\ \{x,y\}\in \mathcal{T}^{\omega_1},\ \forall z\sim x, z\not = y,
 \frac{\omega_1(x,z)}{1-\omega_1(x,y)}=\frac{\omega_2(x,z)}{1-\omega_2(x,y)}.
\]
\end{itemize}
We will denote by $\tilde{\Omega}$ the set of all equivalence classes for the trap-equivalence relation.

\paragraph{Partially forgotten walk.} Set $\tilde{\omega}\in\tilde{\Omega}$. Let $\mathcal{T}$ be its set of traps and $\sigma$ a path starting at $0$ that only stays a finite amount of time every time it enters a trap. We want to define a path whose trajectory is the same as $\sigma$ outside the traps, but in which are erased all the back-and-forths inside traps. For this purpose we define the following sequences of times $(t_i),(s_i)$:
\begin{align*}
t_0&=0,\\
s_i&=\inf\{n\geq t_i, (\sigma_{n}=\sigma_{t_i}\text{ or }\{\sigma_{n},\sigma_{t_i}\}\in\mathcal{T}) \\
&\qquad\qquad\qquad \text{ and } (\sigma_{n+1}\not =\sigma_{t_i}\text{ and }\{\sigma_{n+1},\sigma_{t_i}\}\not\in\mathcal{T})\},\\
t_{i+1}&=
\begin{cases}
s_i+1 \text{ if }\sigma_{s_i}=\sigma_{t_i}\\
s_i \text{ otherwise. }
\end{cases}
\end{align*}
The partially forgotten path $\tilde{\sigma}$ associated with $\sigma$ in the environment $\tilde{\omega}$ is defined by:
\[
\tilde{\sigma}_i:=\sigma_{t_i}.
\]
Similarly we can define the partially-forgotten walk $(\tilde{X}_n)_{n\pg 0}$ associated with $(X_n)_{n\pg 0}$.

Let $f = \{x, y\}$ be an undirected edge of~$\cG$. We first denote by $j$ its direction, that is to say the unique index $1 \pp j \pp d$ such that $y-x = \pm e_j$. Moreover, since the walk is almost surely transient, almost surely the number of visits to~$f$ is finite; let us so define, for every $(u, v) \in \{x, y\}^2$, $N_{u\rightarrow v}$ the number of times the walk enters in $f$ by~$u$ and exits it by~$v$. Therefore, for every unoriented edge $f = \{x, y\}$, we define its configuration $\mathfrak{c}(f)$ by:
$$\mathfrak{c}(f) = (j, N_{x\rightarrow x},N_{x\rightarrow y},N_{y\rightarrow x},N_{y\rightarrow y}).$$
(The choice of which vertex is $x$ or $y$ is not significant, it only needs to be consistent. By convention, we can always assume that $x$ is the vertex at which the walk is located when it hits $f$.)

Note that all the coordinates of $\mathfrak{c}$ only depend on the partially forgotten path, unlike $p_f$, $s_f$ and the $\ell_f^j$ (defined in section~\ref{traps}). Hence, the configuration $\mathfrak{c}(f)$ is a random variable measurable with respect to $(\tilde\omega, \tilde Y)$.

To simplify notations we will also write:
$$N'_x:=N_{x\rightarrow x}+N_{y\rightarrow x} \text{\quad and \quad} N'_y:=N_{x\rightarrow y}+N_{y\rightarrow y}.$$

In \cite{poudevigne2019limit}, Poudevigne has computed the distribution of the strength of the traps conditioned to the partially forgotten walk.

\begin{lemma}[Lemma 2.3.1 of \cite{poudevigne2019limit}]
\label{density_forgotten}
In $\ZZ^d$ endowed with the Dirichlet environment with parameters $(\alpha(e))$, knowing $\tilde{\omega}$ and $(\tilde{X}_i)_{i \pg 0}$, the strength of the various traps are independent.

Moreover, knowing $\tilde{\omega}$ and $(\tilde{X}_i)_{i \pg 0}$, fix a trap $f = \{x,y\}$, and denote by $$\mathfrak{c}=(j, N_{x\rightarrow x},N_{x\rightarrow y},N_{y\rightarrow x},N_{y\rightarrow y})$$ its configuration. Do the change of variables $(\omega(x,y), \omega(y,x)) \to (r,k)$ defined by:
$$(1-\omega(x,y),1-\omega(y,x))=((1+k)r,(1-k)r).$$
The density of law of $(r,k)$ knowing $\tilde{\omega}$ and $\tilde{X}$ is:
\begin{multline*}
C_\mathfrak{c} r^{\kappa_j-1} (1+k)^{N'_x+\alpha_0-\alpha(x,y)-1}(1-k)^{N'_y+\alpha_0-\alpha(y,x)-1} \\*
\times h_\mathfrak{c}(r(1+k),r(1-k))
\Ind{0\leq r\leq\frac{1}{4}}\Ind{-1\leq k \leq 1},
\end{multline*}
where $C_\mathfrak{c}$ is a constant that only depends on the configuration $\mathfrak{c}$, and $h_\mathfrak{c}$ is a function that only depends on $\mathfrak{c}$ and that satisfies the following bound:
\[
\forall 0\pp r\pp \frac{1}{4},\ \ 
\left\vert \log(h_\mathfrak{c}(r(1+k),r(1-k))) \right\vert
\pp 5(N_f+2\alpha_0)r.
\]
\end{lemma}
This lemma can be interpreted as a result of quasi-independence between the strength of a trap (we have $s_f = \frac{1}{2r}$) and its asymmetry~$k$, meaning the imbalance between $\omega(x,y)$ and $\omega(y,x)$. The coupling between $r$ and $k$ is enclosed in the function $h_{\mathfrak{c}}$. The estimates on $h_{\mathfrak{c}}$ show that the intensity of the coupling is essentially controlled by the number of visits to the trap; in other words, to have a decent approximation in saying that $r$ and $k$ are independent, $N_f$ has to be small. This motivates the separation we will make between traps visited a few times and traps visited often.

We deduce easily the following corollary:

\begin{lemma}
\label{quasi-independence}
Keep the notations of the previous lemma.

(i) There also exists a constant $D_\mathfrak{c}$ that only depends on $\mathfrak{c}$ such that for any $A\pg 2$: 
\[
D_\mathfrak{c} A^{-\kappa_j} \exp\left( -\frac{5(N_f + 2\alpha_0)}{2A}\right)
\pp \PP_0\left(s \pg A \big\vert \tilde{\omega},\tilde{Y}\right) 
\pp D_\mathfrak{c} A^{-\kappa_j} \exp\left( \frac{5(N_f + 2\alpha_0)}{2A}\right).
\]
Moreover, there exists a constant $C < \infty$ such that, for all possible configurations~$\mathfrak{c}$,
$$D_\mathfrak{c} \pp C (N_f)^{\kappa_j}.$$

(ii) For any positive concave non-decreasing function $\phi$, denoting $\Phi: t \longmapsto \int\limits_{x=0}^t \phi(x) \dd x$, there exists a constant $C_\phi$, such that:
\[
\forall A\pg 2, \qquad \sE_0\left[\Phi(N_f) \Ind{N_f \pg 1} \Ind{s_f\pg A}\right] \pp C_\phi \cdot \sE_0\left[\Phi(N_f)\right] {A^{-\kappa_j}}
\]
for every unoriented edge $f$ colinear to $e_j$.
\end{lemma}

\begin{proof} We only prove (i), since (ii) is already proved in of \cite{poudevigne2019limit} (lemma~2.3.4).

(i) is a simple consequence of the previous lemma~\ref{density_forgotten}. Set:
$$ D_\mathfrak{c} = \frac{1}{\kappa_j 2^{\kappa_j}} \int_{-1}^1 C_\mathfrak{c} (1+k)^{N'_x+\alpha_0-\alpha(x,y)-1}(1-k)^{N'_y+\alpha_0-\alpha(y,x)-1} \dd k.$$

By lemma~\ref{density_forgotten}:
\[
\forall 0\pp r\pp \frac{1}{4},
e^{-5(N_f+2\alpha_0)r} \pp 
h_\mathfrak{c}(r(1+k),r(1-k))
\pp e^{5(N_f+2\alpha_0)r},
\]
so, by integrating the density of the law of $(r,k)$ over $\{s \pg A\}$, i.e. over $\left\lbrace r \pp \frac{1}{2A} \right\rbrace$:
\[
D_\mathfrak{c} A^{-\kappa_j} \exp\left( -\frac{5(N_f + 2\alpha_0)}{2A}\right)
\pp \PP_0\left(s \pg A \big\vert \tilde{\omega},\tilde{Y}\right) 
\pp D_\mathfrak{c} A^{-\kappa_j} \exp\left( \frac{5(N_f + 2\alpha_0)}{2A}\right).
\]

But we know that:
$$ 1 = \int_{r=0}^{1/4} \int_{k=-1}^1 C_\mathfrak{c} r^{\kappa_j-1} (1+k)^{N'_x+\alpha_0-\alpha(x,y)-1}(1-k)^{N'_y+\alpha_0-\alpha(y,x)-1} h_{\mathfrak{c}}(r(1-k), r(1+k)) \dd k \dd r$$
with, for every $r > 0$:
$$ h_\mathfrak{c}(r(1+k),r(1-k)))
\pg e^{-5(N_f+2\alpha_0)r}; $$
so that:
\begin{align*}
1 &\pg \int_{r=0}^{1/4} \left( \int_{k=-1}^1 C_\mathfrak{c} (1+k)^{N'_x+\alpha_0-\alpha(x,y)-1}(1-k)^{N'_y+\alpha_0-\alpha(y,x)-1} \dd k \right) e^{-5(N_f+2\alpha_0)r} r^{\kappa_j-1} \dd r \\
&\pg \kappa_j 2^{\kappa_j} D_\mathfrak{c} \cdot \int_{r=0}^{1/4} e^{-5(N_f+2\alpha_0)r} r^{\kappa_j-1} \dd r \\
&\pg \kappa_j 2^{\kappa_j} D_\mathfrak{c}  \left(5(N_f+2\alpha_0)\right)^{-\kappa_j} \int_{r=0}^{\frac{5}{4}(N_f + 2\alpha_0)} e^{-r} r^{\kappa_j-1} \dd r \\
&\pg D_\mathfrak{c} \cdot C N_f^{-\kappa_j}
\end{align*}
This concludes the proof.
\end{proof}

The main consequence of the previous lemma is the following estimate of the tail of the time spent in a given trap, proved in~\cite{poudevigne2019limit}:

	\begin{lemma}[proved in the proof of lemma 2.4.3 of \cite{poudevigne2019limit}]
	\label{tail}
	Let $f$ be an unoriented edge. Let $\mathfrak{c} = (j, N_{x\rightarrow x}, N_{x\rightarrow y}, N_{y\rightarrow y}, N_{y\rightarrow x})$ be its configuration. Fix some function $h: \RR_+ \longmapsto \RR_+$ which tends to $\infty$ at $+\infty$, but at most linearly. There exists a constant $C_\mathfrak{c}$ such that, when $x \to \infty$ we have the equivalent:
	$$ \PP\left(\left.\sum_{i=1}^{N_f} \ell_f^i \pg x, s_f \pg h(x) \right\vert \tilde\omega, \tilde X \right) \sim C_{\mathfrak{c}} x^{-\kappa_j}. $$
	\end{lemma}
	
	\subsection{The time spent out of traps}
	
We begin by discarding the contribution of the time spent out of traps in the tail of $T_1$, which we expect to be of order $x^{-\kappa}$. We first show that estimating $\sE_x[(T_1^\circ)^\beta]$ reduces to estimating the order of integrability of Green functions in a new environment. We then prove in the following lemma that the Green functions have higher order of integrability in this new environment that in the original one, enabling us to take $\beta > \kappa$ and keep $\sE_x[(T_1^\circ)^\beta]$ finite.

\begin{lemma}
\label{out_lemma1}
For every $\beta > 1$, there exists a constant $C_\beta$ such that, for every $r>1$ and $s>0$:
\begin{multline*}
\sE_x\left[(T_1^\circ)^\beta \Ind{T_1 < \bar H_U} \right]
\pp (8d-4)^{r} (2d)^s C_\beta |U|^{\beta-1} \sum_{\substack{y \in U, \underline{e}=y, \\ I \in \cI_y}} \EE^{(\alpha)}\left[ \omega(e)^s (P^U_y)^{-\beta} \prod_{e' \in I} \omega(e')^{r/2d}\right]
\end{multline*}
where $P_y^U = \sP_y^\omega(H^+_y > \bar H_U)$ and $\cI_y$ is a set of $2d$-tuples of oriented edges which has the following properties: for every $I \in \cI_y$ and $x \sim y$, there exists $e' \in I$ such that $e' \in \partial_+\{x,y\}$.
\end{lemma}

\begin{proof}
Let $\beta > 0$. For every $r > 0$, by Markov's inequality, we have for every $\omega \in \Omega$:
\begin{align*}
\sE_x^\omega\left[(T_1^\circ)^\beta\Ind{T_1 < \bar H_U} \right]
&= \sE_x^\omega\left[\left(\sum_{i=1}^{T_1} \Ind{X_{i-1} \in U\setminus\tilde{\cT}}\right)^\beta\right] \\
&= \sE_x^\omega\left[\left(\sum_{i=1}^{T_1} \sum_{y \in U} \Ind{X_{i-1} = y} \Ind{y \not\in \tilde{\cT}} \right)^\beta\right] \\
&\pp |U|^{\beta-1} \sum_{y \in U} \sE_x^\omega\left[\left(\sum_{i=1}^{T_1} \Ind{X_{i-1} = y} \right)^\beta \Ind{s_y \pp 2}\right] \\
&\pp 2^r |U|^{\beta-1} \sum_{y \in U} s_y^{-r} \sE_x^\omega\left[\left(\sum_{i=1}^{T_1} \Ind{X_{i-1} = y} \right)^\beta \right].
\end{align*}

By Markov property:
\begin{align*}
\sE_x^\omega\left[(T_1^\circ)^\beta\Ind{T_1 < \bar H_U} \right]
&\pp 2^r |U|^{\beta-1} \sum_{y \in U} s_y^{-r} \sP_x^\omega(H_y < \bar H_U) \sE_y^\omega\left[\left(\sum_{i=1}^{T_1} \Ind{X_{i-1} = y} \right)^\beta \right] \\
&\pp 2^r |U|^{\beta-1} \sum_{y \in U} s_y^{-r} \sE_y^\omega\left[\left(\sum_{i=1}^{T_1} \Ind{X_{i-1} = y} \right)^\beta \right].
\end{align*}

Note that, under $\sP^\omega$, for every $y \in U$, the number $\sum_{i=1}^{T_1} \Ind{X_{i-1} = y}$ of visits of the walk at~$y$ before it exits the box~$U$ follows a geometric distribution with parameter:
$$ P_y^U = \sP_y^\omega(H^+_y > \bar H_U) ; $$
therefore, by lemma~\ref{appendix}, there exists a finite constant $C_\beta$ such that:
$$ \sE_y^\omega\left[\left(\sum_{i=1}^{T_1} \Ind{X_{i-1} = y} \right)^\beta \right] \pp C_\beta (P^U_y)^{-\beta}. $$

Integrating the last inequality over $\omega \in \Omega$ gives finally, for every $r > 0$:
$$ \sE_x\left[(T_1^\circ)^\beta\right]
\pp 2^r C_\beta |U|^{\beta-1} \sum_{y \in U} \EE^{(\alpha)}\left[ s_y^{-r} (P^U_y)^{-\beta}\right]. $$

Recall that, for every $y \in V$,
$$ s_y = \max_{x\sim y} \left( \frac{1}{2-\omega(x,y)-\omega(y,x)}\right) = \left( \min_{x\sim y} \sum_{e \in \partial_+\{x, y \}} \omega(e) \right)^{-1}. $$
We can thus, in the previous inequality, use
$$ s_y^{-1} \pp \prod_{x\sim y} \left(\sum_{e' \in \partial_+\{x, y \}} \omega(e') \right)^{1/2d}. $$
Moreover, using a trick already used in~\cite{sabot2013random}, we formulate mathematically that a single vertex cannot be a trap (because the walk necessarily steps out of this vertex):
$$ \sum_{e\in E \tq \underline{e}=y} \omega(e) = 1.$$

Therefore, for every $s > 0$,
$$ \sE_x\left[(T_1^\circ)^\beta\right]
\pp 2^r C_\beta |U|^{\beta-1} \sum_{y \in U} \EE^{(\alpha)}\left[ \prod_{x\sim y} \left(\sum_{e' \in \partial_+\{x, y \}} \omega(e') \right)^{r/2d} \left(\sum_{\underline{e}=y} \omega(e)\right)^s (P^U_y)^{-\beta}\right]. $$
Since all quantities here are non-negative, for every $s > 0$, a convexity inequality gives:
\begin{multline*}
\sE_x\left[(T_1^\circ)^\beta\right]
\pp 2^r (2d)^s (4d-2)^{r/2d} C_\beta |U|^{\beta-1} \sum_{y \in U} \EE^{(\alpha)}\left[ \prod_{x\sim y} \left( \sum_{e' \in \partial_+\{x, y \}} \omega(e')^{r/2d} \right) \right. \\*
\left. \cdot \sum_{\underline{e}=y} \omega(e)^s (P^U_y)^{-\beta}\right].
\end{multline*}
Let us finally develop the product, and more specifically take caution to the inversion of the indices: given a vertex~$y$, choosing $x$ a neighbour of~$y$ and then choosing one single edge $e' \in \partial_+\ens{x,y}$ is equivalent to choose a collection~$I$ of~$2d$ edges in~$\cI_y$ and then choosing an edge~$e'$ in~$I$. Therefore:
$$\prod_{x\sim y} \left( \sum_{e' \in \partial_+\{x, y \}} \omega(e')^{r/2d} \right) = \sum_{I \in \cI_y} \prod_{e'\in I} \omega(e')^{r/2d} .$$

Finally,
$$ \sE_x\left[(T_1^\circ)^\beta\right]
\pp (2d)^s (8d-4)^{r} C_\beta |U|^{\beta-1} \sum_{y \in U} \sum_{\underline{e}=y} \sum_{I \in \cI_y} \EE^{(\alpha)}\left[\prod_{e'\in I} \omega(e')^{r/2d} \omega(e)^s (P^U_y)^{-\beta}\right]. $$
\end{proof}

\begin{lemma}
\label{out_lemma2}
There exists $\beta > \kappa$ and a constant $C_\beta < \infty$ such that:
$$ \sE_0\left[(T_1^\circ)^\beta \Ind{T_1 < \bar H_U} \right] < C_\beta \vert U \vert^{\beta}. $$
\end{lemma}

\begin{proof}
Let $\beta > 1$. Fix $y \in U$, $e\in E$ such that $\underline{e} = y$ and $I = (e_1, \ldots, e_{2d}) \in \cI_y$. The strategy of the proof is the following. Sabot and Tournier showed in~\cite{sabot2011reversed} the following consequence of the statistical invariance by time reversal of the Dirichlet environment on finite graphs with null divergence:

Consider the graph $\cG_U = (U \cup \{\partial\}, E_U)$ obtained from $\ZZ^d$ by contracting all the points outside from~$U$ in a single point~$\partial$ and adding an edge with weight $\gamma$ from $\partial$ to $y$; choose $\theta$ a unitary flow from $y$ to $\partial$ so that $\alpha+\gamma\theta$ has null divergence; then, under $\PP^{(\alpha + \gamma\theta)}$, $P_y^U$ is stochastically lower-bounded by a beta distribution $\text{Beta}(\gamma, \alpha_y)$.

We therefore perform this change of measure in the bound of lemma~\ref{out_lemma1} (with a different $\theta = \theta_{y, e, I}$ for each expectation). On the one hand, we must check that, under this new environment $\PP^{(\alpha+\gamma\theta)}$, the expectations in the bound of lemma~\ref{out_lemma1} do not explode, thanks to the max-flow min-cut theorem. On the other hand, we must check that the cost of these changes of measure is the same for every expectation: this is ensured by global estimates of the $L^2$-norm of unitary flows in $\ZZ^d, d \pg 3$. 

\paragraph{Change of measure.} We thus proceed to the change of measure from $\PP^{(\alpha)}$ to $\PP^{(\alpha + \gamma\theta)}$: by equation~\ref{eq:changement_mesure}, for every $r>1$ and $s>0$,
$$ \EE^{(\alpha)}\left[\prod_{e'\in I} \omega(e')^{r/2d} \omega(e)^s (P^U_y)^{-\beta}\right] = \frac{Z_{a+\gamma\theta}}{Z_\alpha} \EE^{(\alpha+\gamma\theta)}\left[\prod_{e'\in I} \omega(e')^{r/2d} \omega(e)^s (P^U_y)^{-\beta} \omega^{-\gamma\theta}\right]. $$

By Hölder's inequality, for every $p, q > 1$ such that $\frac{1}{p} + \frac{1}{q} = 1$, we thus have:
\begin{align*}
\EE^{(\alpha)}\left[\prod_{e'\in I} \omega(e')^{r/2d} \omega(e)^s (P^U_y)^{-\beta}\right]
\pp {} & \overbrace{\EE^{(\alpha+\gamma\theta)}\left[ (P^U_y)^{-p\beta} \right]^{1/p}}^{\boxed 1} \\
& {}\cdot \underbrace{\frac{Z_{a+\gamma\theta}}{Z_\alpha} \, \EE^{(\alpha+\gamma\theta)}\left[ \prod_{e'\in I} \omega(e')^{qr/2d} \omega(e)^{qs} \omega^{-\gamma q\theta} \right]^{1/q}}_{\boxed 2}
\end{align*}

Remark that we have introduced four independent parameters $q > 1$, $r > 0$, $s > 0$ and $\gamma > 0$. Until now, these parameters need not to satisfy other conditions.

Let us examine each line of the last bound.

\paragraph{First line: statistical invariance by time reversal.} Applying the result of~\cite{sabot2011reversed} reminded above, we have:
$$ \boxed{1} \pp \EE\left[ \text{Beta}(\gamma, \alpha_y)^{-p\beta} \right]. $$
It suffices that $p\beta < \gamma$ for this upper-bound, which is independent on~$U$ and~$y$, to be a finite constant. This is the first condition between our parameters.

\paragraph{Second line: max-flow min-cut theorem.} The expectation
$$ \EE^{(\alpha+\gamma\theta)}\left[ \omega(e)^{qs} \left(\prod_{e'\in I} \omega(e')^{qr/2d}\right) \omega^{-\gamma q\theta} \right] $$
is finite if and only the unitary flow $\theta = \theta_{y, e, I}$ is such that:
$$ \Div(\theta) \equiv \delta_y - \delta_\partial $$
and such that, for every oriented edge $f \in E_U$,
$$ \theta(f) \pp \frac{1}{\gamma(q-1)} \left[\alpha(e) + qs \Ind{e}(f) + \frac{qr}{2d} \sum_{e'\in I} \Ind{e'}(f) \right]. $$
By an extension of the max-flow min-cut theorem, these two conditions are compatible as soon as we can find some conductances $(c(f)) = (c_{y, e, I}(f))_{f\in E_U}$ such that:
\begin{enumerate}
\item $m(c):= \inf \left\lbrace \sum_{f \in \partial_+ A} c(f), A \text{ connected subset of } U \text{ containing }y \right\rbrace \pg 1$;
\item for every $f \in E_U$,
$$ c(f) \pp \frac{1}{\gamma(q-1)} \left[\alpha(f) + qs \Ind{e}(f) + \frac{qr}{2d} \sum_{e' \in I} \Ind{e'}(f) \right]. $$
\end{enumerate}

Let us prove these two properties are verified by the following definition of $c$:
$$\forall f \in E_U, \quad c(f) = \delta \left[\alpha(f) + qs \Ind{e}(f) + \frac{qr}{2d} \sum_{e' \in I} \Ind{e'}(f) \right]$$
where $\delta > 0$ is to be chosen judiciously later.

Let $A$ be a subset of $U$ connected at~$y$ which realizes the minimum in $m(c)$. It is one of two things:
\begin{itemize}
\item either $\partial_+ A$ contains an edge among $I \cup \{e\}$; in this case,
$$m(c) \pg \inf_{f \in I \cup \{e\}} c(f) \pg \delta \left(\min \{\alpha_1, \ldots, \alpha_{2d}\}  + \min\left\lbrace qs, \frac{qr}{2d} \right\rbrace \right); $$
\item or $\partial_+ A$ contains none of these edges; since $A$ is connected at~$y$, and $e \not\in \partial_+ A$, then $\bar{e} \in A$; moreover, $\bar{e} \sim y$ so, by construction of $I$, there exists an $e'\in I$ such that $e'$ goes out from the set $\{y, \bar{e}\}$, and, by the same reasoning, $\bar{e'} \in A$---finally, $A$ contains three distinct points $\{x, y, z\}$ such that $x \sim y \sim z$, so
$$ m(c) \pg \delta \sum_{e \in \partial_+ A} \alpha(e) \pg \delta \kappa(\{x, y, z\}) \pg \epsilon \kappa'. $$
\end{itemize}

Therefore, we must choose
$$ \frac{1}{\kappa'} < \delta < \frac{1}{\gamma (q-1)}, \qquad s > \frac{1}{q\delta}, \qquad r > \frac{2d}{q\delta}. $$
in order $c$ to satisfy the two conditions aforementioned.

We finally show that all these inequalities are compatible.

\paragraph{Recapitulation: final choice of parameters.} Finally, we must choose all the parameters we introduced in the proof as follows:
\begin{enumerate}
\item let $q \in \left(1, \frac{\kappa'}{\kappa}\right)$ and $p$ be its conjugated exponent;
\item let $\gamma \in \left( \frac{q}{q-1} \kappa, \frac{\kappa'}{q-1}\right)$---which is possible thanks to the previous condition;
\item let $\delta \in \left(\frac{1}{\kappa'}, \frac{1}{\gamma (q-1)}\right)$---which is possible thanks to the previous condition;
\item let $s > \frac{1}{q\delta}$ and $r > \frac{2d}{q\delta}$.
\end{enumerate}
With this choice of parameters (independent on~$U$ and~$y$), for every $\beta \in \left(\kappa, \frac{\gamma}{p}\right)$---this interval being non-empty,--- $\boxed{2}$ is finite.

\paragraph{Cost of the change of measure: global $L^2$-norm of unitary flows.} Fix $\beta \in \left(\kappa, \frac{\gamma}{p}\right)$.

We already know that $\boxed{1}$ is upper-bounded by a finite constant independent on~$U$ and~$y$.

On the other hand, we have for every $y \in U$, $\underline{e} = y$ and $I \in \cI_y$:
\begin{align*}
\boxed{2} &= \frac{\left(Z_{\alpha - \gamma(q-1)\theta_{y, e, I} + qs \Ind{e} + \frac{qr}{2d} \sum_{e'\in I} \Ind{e'}}\right)^\frac{1}{q}}{Z_\alpha \cdot (Z_{\alpha+\gamma\theta_{y, e, I}})^{\frac{1}{q}-1}} \\
&= \frac{\prod_{x \in U} \Gamma(\alpha(x))}{\prod_{f \in E_U} \Gamma(\alpha(f))} \cdot \left(\frac{\prod_{x \in U} \Gamma(\alpha(x) + \gamma\theta_{y, e, I}(x))}{\prod_{f \in E_U} \Gamma(\alpha(f) + \gamma\theta_{y, e, I}(f)))}\right)^{\frac{1}{q}-1} \\
&\qquad \cdot \left(\frac{\prod_{f\in E_U} \Gamma\left(\alpha(f) - \gamma(q-1)\theta_{y, e, I}(f) + qs \Ind{e}(f) + \frac{qr}{2d} \sum_{e'\in I} \Ind{e'}(f)\right)}{\prod_{x\in U} \Gamma\left(\alpha(x) - \gamma(q-1)\theta_{y, e, I}(x) + qs \Ind{y}(x) + \frac{qr}{2d} \sum_{e'\in I} \Ind{\underline{e'}}(x)\right)} \right)^{\frac{1}{q}}
\end{align*}
By construction of $\theta_{y, e, I}$, for every edge $f \in I_e:= \{e\} \cup I$ and every $x \in \underline{I}_e:= \ens{\underline{e'}, e'\in I_e}$,
\begin{align*}
\Gamma(\alpha(f)) &\pg \inf\ens{\Gamma(s), s\in [\min_i \alpha_i, \max_i \alpha_i]}  \\
\Gamma(\alpha(x)) &= \Gamma(\alpha_0) \\
\Gamma\left(\alpha(f) + \gamma\theta_{y, e, I}(f)\right) &\pg \inf\ens{\Gamma(s), s\in [\min_i \alpha_i, (1+p) \max_i \alpha_i + pr + ps]} \\
\Gamma\left(\alpha(x) + \gamma\theta_{y, e, I}(x)\right) &\pp \sup\ens{\Gamma(s), s\in [\alpha_0, (1+p) \alpha_0 + pr + ps]}
\end{align*}
\begin{multline*}
\Gamma\left(\alpha(f) - \gamma(q-1)\theta_{y, e, I}(f) + qs \Ind{e}(f) + \frac{qr}{2d} \sum_{e' \in I} \Ind{e'}(f)\right) \\
\pp \sup\ens{\Gamma(s), s\in [(1-\gamma(q-1)\delta) \min_i \alpha_i, \max_i \alpha_i + qr + qs]}
\end{multline*}
\begin{multline*}
\Gamma\left(\alpha(x) - \gamma(q-1)\theta_{y, e, I}(x) + qs \Ind{y}(x) + \frac{qr}{2d} \sum_{e' \in I} \Ind{\underline{e_i}}(x)\right) \\
\pg \inf\ens{\Gamma(s), s\in [(1-\gamma(q-1)\delta) \min_i \alpha_i, \alpha_0 + qr + qs]}
\end{multline*}
All these bound are strictly positive, finite and independent from~$U$ and~$y$. Therefore, isolating the contributions of the edges $f \in I_e$ and the vertices $x \in \underline I_e$ in $\boxed{2}$, there exists a finite constant $A$, depending only on $d$, $(\alpha_1, \ldots, \alpha)$ and $(q, \gamma, \delta, r, s)$ and thus independent on $U$ and $y$ such that:
\begin{align*}
\boxed{2} &\pp A \cdot \frac{\prod_{x \not\in \underline I_e} \Gamma(\alpha(x))}{\prod_{f \not\in I_e} \Gamma(\alpha(f))} \cdot \left(\frac{\prod_{x \not\in \underline I_e} \Gamma(\alpha(x) + \gamma\theta_{y, e, I}(x))}{\prod_{f \not\in I_e} \Gamma(\alpha(f) + \gamma\theta_{y, e, I}(f)))}\right)^{\frac{1}{q}-1} \\
&\qquad \cdot \left(\frac{\prod_{f\not\in  I_e} \Gamma\left(\alpha(f) - \gamma(q-1)\theta_{y, e, I}(f)\right)}{\prod_{x\not\in \underline I_e} \Gamma\left(\alpha(x) - \gamma(q-1)\theta_{y, e, I}(x)\right)} \right)^{\frac{1}{q}} \\
&\pp A \exp\left(\sum_{f \not\in I_e} \nu_{q}(\alpha(f), \gamma\theta(f)) - \sum_{x \not\in \underline I_e} \nu_{q}(\alpha(x), \gamma\theta(x)) \right)
\end{align*}
where for every $\alpha > 0$ and $t < \frac{\alpha}{1-q}$:
\begin{align*}
\nu_{q}(\alpha, t) = \frac{1}{q} \ln\Gamma(\alpha - (q-1)t) + \left(1-\frac{1}{q}\right) \ln\Gamma(\alpha + t) -  \ln\Gamma(\alpha).
\end{align*}
By Taylor's inequality, there exists a constant $a$ such that, for every $\alpha$ in a fixed compact set of $\RR_+^*$ containing $\alpha(0), \alpha_1, \ldots, \alpha_{2d}$ and $t \pp \epsilon\alpha$:
$$ \vert \nu_{\gamma, q}(\alpha, t) \vert \pp a \vert t \vert^2. $$
Since we have chosen all $\theta_{y, e, I}$ so that $\theta_{y, e, I}(f) \pp c_{y, f, (e_i)}(f) \pp \delta\alpha(f)$ for every $f \not\in I_e$, this inequality gives:
$$ \boxed{2} \pp A \exp\left(a\left[ \sum_{e\in E_U} \theta_{y, e, I}^2 + \sum_{x \in U} \theta_{y, e, I}^2(x) \right]\right). $$

Since our conductances $(c_{y, f, (e_i)}(e))$ are bounded:
$$\frac{1}{\gamma(q-1)} \min_i\alpha_i \pp c_{y, f, (e_i)}(e) \pp \frac{1}{\gamma(q-1)} (\max_i\alpha_i + qr + qs),$$
and since these bounds are independent on $U$ and $y$, the same lemma we used earlier ensures that the $L^2$-norms of all $\theta_{y, e, I}$ are of same order of the electrical resistance between 0 and~$\ZZ\setminus U$ for the network~$\ZZ^d$ with unit resistance on the bonds---in particular, there exists a finite constant $a'$ such that these norms are bounded by~$a'$ when $d\pg 3$: so there exists finite constants $A'$ and $\beta_0$ independent of $U$ and $y$ such that
$$ \boxed{2} \pp A' |U|^{\beta_0}. $$

Finally there exists a finite constant $A''$ such that, for every $y\in U$, $\underline e = y$ and $I \in \cI_y$,
$$ \EE^{(\alpha)}\left[\prod_{e'\in I} \omega(e')^{r/2d} \omega(e)^s (P^U_y)^{-\beta}\right] \pp A'' |U|^{\beta_0}. $$
This concludes the proof thanks to the previous lemma.
\end{proof}

\subparagraph*{Remark 1.} If we have the same reasoning for the random walk in Dirichlet environment in dimension ${d = 2}$, we obtain the following bound: there exists $\beta > \kappa$, $\beta' > 0$ and a constant $C_\beta < \infty$ such that:
$$ \sE_0\left[(T_1^\circ)^\beta \Ind{T_1 < \bar H_U} \right] < C_\beta \vert U \vert^{\beta+\beta'}. $$
The extra power appears because, in dimension~2, the norms of the flows $\theta_{y,e,I}$ are not uniformly bounded by a constant, but by $a'\ln|U|$ for some $a'>0$ (cf. \cite{sabot2017random}).
	
	\subsection{The number of visits of traps}
	
	We turn now to the estimation of the number of visits to a given trap. We proceed by contracting this trap in a single point.

	\paragraph{Contraction lemma.} The effect of contraction on Dirichlet distribution has been studied in~\cite{tournier2009integrability}:
	
	\begin{lemma}
	\label{contraction}
Let $r$, $(d_1, \ldots, d_r)$ and $(\delta_1, \ldots, \delta_r)$ be positive integers such that $\delta_i < d_i$ for every $1 \pp i \pp r$, and $(\alpha_{x,y}, 1\pp x \pp r, 1 \pp y \pp d_x)$ some positive real numbers. Assume that, independently for every~$1 \pp x \pp r$, the vector $(\omega_{x,y})_{1 \pp y \pp d_x}$ follows $\cD(\alpha_{x,y}, 1 \pp y \pp d_x)$.

Denote $\beta = \sum_{x=1}^r \sum_{y=1}^{\delta_x} \alpha_{x,y}$, $\Sigma = \sum_{x=1}^r \sum_{y=1}^{\delta_x} \omega_{x,y}$ and $\Delta = \sum_{x=1}^r \delta_x$.

There exists two finite constants $C, C' > 0$ such that, for every function $f: \RR\times\RR^\Delta \to \RR$, 
$$E\left[f\left(\Sigma, \left(\frac{\omega_{x,y}}{\Sigma}\right)_{\substack{1\pp x \pp r, \\ 1 \pp y \pp \delta_x}}\right)\right]
	\pp C \cdot \tilde{E}\left[f\left({\Sigma}^{\mathrm{Dir}},(\omega_{x,y}^{\mathrm{Dir}})_{\substack{1\pp x \pp r, \\ 1 \pp y \pp \delta_x}}\right)\right],$$
where, under the probability $\tilde{P}$, $\tilde{\Sigma}$ is a random variable independent from the $\omega^{\mathrm{Dir}}$ and the random vector $(\omega_{x,y}^{\mathrm{Dir}})_{\substack{1\pp x \pp r, \\ 1 \pp y \pp \delta_x}}$ follows $\cD(\alpha_{x,y}, 1\pp x \pp r, 1 \pp y \pp \delta_x)$. 
\end{lemma}

What is the sense of lemma~\ref{contraction}? Consider a graph $\cG = (V, E)$ endowed with a Dirichlet environment with weights $(\alpha(e))_{e\in E}$. Let $f \subset V$ be a connected set of vertices. Denote $r = |f|$ and call conventionnally $1, \ldots, r$ the different vertices of~$f$. Contract $f$ in a single point, \emph{i.e.} construct a new graph $\cG^f = (V^f, E^f)$ where $V^f$ is the reunion of $V\setminus f$ and a new vertex $x_f$, and where $E^f$ is composed of the edges of $E$ between vertices of $V\setminus f$, to which we add an edge between a vertex of $v \in V\setminus f$ and $x_f$ (respectively $x_f$ and $v \in V\setminus f$) every time that there exists a vertex $u \in f$ such that $(v,u)\in E$ (respectively $(u, v)\in E$)---this can lead to the creation of multiple edges between two vertices. We denote by $d_i$ the numbers of neighbours of vertex~$i$ in $\cG$, and $\delta_i$ the number of neighbours of vertex~$i$ in $\cG$ that are not elements of~$f$. $\Delta$ is then the degree of $x_f$ in $\cG^f$.

For every environment $\omega$ over $\cG$, we define an environment $\tilde\omega$ over $\cG^f$. We conserve $\omega(e)$ when $e$ is an edge such that $\underline{e} \in V \setminus f$, but, on edges from $x_f$, corresponding to an original edge from a vertex~$i$ of $f$ to a vertex $v \in V\setminus f$, we define $\tilde\omega(e)$ as $\omega(i, v)$ renormalized by $\Sigma$.

Because of the renormalization factor, this new environment is not a Dirichlet environment over~$\cG^f$. This lemma, however, ensures that it can be compared to a Dirichlet environment over~$\cG^f$.

\paragraph{Number of visits in a trap.} We use this contraction lemma to compare the number of visits in a trap and the Green's function in the graph where this point has been contracted in a single point.
	
\begin{lemma}
\label{visits1}
In $\ZZ^d, d \pg 3$, there exists $\beta > \kappa$ and a constant $C_\beta \in (0,\infty)$, such that, for every $f = \{x, y\}$ an unoriented edge:
$$ \sE_0\left[(N_f)^\beta \right] \pp C_\beta. $$
\end{lemma}

\begin{proof}
Remind that if $f$ is an unoriented edge, we have denoted by $N_f$ the number of visits of the walk to $f$ before time~$T_1$. If $x$ is a vertex, we may denote by $N_x$ the number of visits of the walk at $x$ before time~$T_1$. In a similar way, we now denote by $N_{f \vert U}$ and $N_{x \vert U}$ the numbers of visits of the walk to $f$ or $x$ before it exits a given subset of vertices~$U$.

Let $U$ be any given subset of $V$ containing $f$. It is a classical result that, since $X$ is a Markov chain under the probability~$\sP^\omega$, then under~$\sP_x^\omega$, $N_{x \vert U}$ has a geometric distribution with parameter $\sP^\omega(H_x > \bar{H}_U)$. To control $\sE_0\left[(N_{f \vert U})^\beta \right]$, we aim at a similar result for $N_{f \vert U}$.

First, it is sufficient to control $\EE\left[ \sup_{v \in f} \sE^\omega_v\left[(N_{f \vert U})^\beta\right] \right]$. Indeed, by strong Markov's property (under~$\sP_x^\omega$), we have, for every $\beta > 1$:
\begin{align*}
\sE_0\left[(N_{f \vert U})^\beta\right]
&= \sE_0\left[\Ind{H_f < \bar H_U} (N_{f \vert U})^\beta \right] \\
&= \EE\left[\sE_0^\omega \left[ \Ind{H_f < \bar H_U} \sE^\omega_{X_{H_f}}\left[(N_{f \vert U})^\beta\right] \right] \right] \\
&= \sum_{u \in f} \EE\left[\sE_0^\omega \left[ \Ind{H_f < \bar H_U, X_{H_f} = u} \sE^\omega_u\left[(N_{f \vert U})^\beta\right] \right] \right] \\
&= \sum_{u \in f} \EE\left[\sP_0^\omega \left(H_f < \bar H_U, X_{H_f} = u \right) \sE^\omega_u\left[(N_{f \vert U})^\beta\right] \right] \\
&\pp \EE\left[\left( \sum_{u \in f} \sP_0^\omega \left(H_f < \bar H_U, X_{H_f} = u \right) \right) \sup_{v \in f} \sE^\omega_v\left[(N_{f \vert U})^\beta\right] \right] \\
&\pp \EE\left[\sP_0^\omega \left(H_f < \bar H_U\right) \sup_{v \in f} \sE^\omega_v\left[(N_{f \vert U})^\beta\right] \right] \\*
&\pp \EE\left[ \sup_{v \in f} \sE^\omega_v\left[(N_{f \vert U})^\beta\right] \right]
\end{align*}

We now reduce the problem to the control of the $\sup_{v\in f} \left(\sP_v^\omega(N_{f \vert U} \pg k)\right)_{k \pg 1}$. A standard application of Fubini's theorem shows that, for $v \in f$, $\beta > 1$ and $\omega\in\Omega$:
$$ \sE^\omega_v\left[(N_{f \vert U})^\beta\right]
= \int_0^{+\infty} \beta t^{\beta-1} \sP_v^\omega(N_{f \vert U} > t) \, \dd t
= \sum_{k=0}^{\infty} \int_k^{k+1} \beta t^{\beta-1} \sP_v^\omega(N_{f \vert U} > t) \, \dd t.
$$
Note that, for every $t \in [k, k+1)]$, $\sP_v^\omega(N_{f \vert U} > t) = \sP_v^\omega(N_{f \vert U} \pg k+1)$. Therefore:
$$ \sE^\omega_v\left[(N_{f \vert U})^\beta\right]
= \sum_{k=0}^{\infty} \int_k^{k+1} \beta t^{\beta-1} \sP_v^\omega(N_{f \vert U} \pg k+1) \, \dd t
= \sum_{k=0}^{\infty} \left[(k+1)^\beta - k^\beta\right] \sP_v^\omega(N_{f \vert U} \pg k+1),
$$
and, returning to our expectation of interest:
$$ \sE_0\left[(N_{f \vert U})^\beta\right]
\pp \EE\left[\sup_{v \in f} \sum_{k=0}^{\infty} \left[(k+1)^\beta - k^\beta\right] \sP_v^\omega(N_{f \vert U} \pg k+1) \right].
$$

Let us so determine the growth order of $\left(\sP_v^\omega(N_{f \vert U} \pg k+1)\right)_{k\pg 0}$. For $k \pg 1$, $v \in f$ and $\omega\in\Omega$, by strong Markov's property, we have:
\begin{align*}
\sP_v^\omega(N_{f \vert U} \pg k+1)
&= \sP_v^\omega(H_f^+ < \bar H_U, N_{f \vert U} \pg k+1) \\
&= \sE_v^\omega\left[ \Ind{H_f^+ < \bar H_U} \sP_{X_{H_f^+}}^\omega (N_{f \vert U} \pg k) \right] \\
&= \sum_{u \in f} \sE_v^\omega\left[ \Ind{H_f^+ < \bar H_U, X_{H_f^+} = u} \sP_u^\omega (N_{f \vert U} \pg k) \right] \\
&= \sum_{u \in f} \sP_v^\omega( H_f^+ < \bar H_U, X_{H_f^+} = u) \sP_u^\omega (N_{f \vert U} \pg k) \\
&\pp \left(\sum_{u \in f} \sP_v^\omega( H_f^+ < \bar H_U, X_{H_f^+} = u)\right) \sup_{u \in f} \sP_u^\omega (N_{f \vert U} \pg k) \\
&\pp \sP_v^\omega( H_f^+ < \bar H_U) \cdot \sup_{u \in f} \sP_u^\omega (N_{f \vert U} \pg k) \\
&\pp \sup_{u \in f} \sP_u^\omega( H_f^+ < \bar H_U) \cdot \sup_{u \in f} \sP_u^\omega (N_{f \vert U} \pg k)
\end{align*}
so that finally:
$$ \sup_{v \in f} \sP_v^\omega(N_{f \vert U} \pg k+1) \pp \sup_{u \in f} \sP_u^\omega( H_f^+ < \bar H_U) \cdot \sup_{v \in f} \sP_v^\omega(N_{f \vert U} \pg k). $$
This directly implies that:
$$ \sup_{v \in f} \sP_v^\omega(N_{f \vert U} \pg k+1) \pp \left[ \sup_{v \in f} \sP_v^\omega( H_f^+ < \bar H_U) \right]^k $$
and, returning to our expectation of interest:
$$ \sE_0\left[(N_{f \vert U})^\beta\right]
\pp \EE\left[\sum_{k=0}^{\infty} \left[(k+1)^\beta - k^\beta\right] \left[ \sup_{v \in f} \sP_v^\omega( H_f^+ < \bar H_U) \right]^k \right].
$$
We are so led to estimating $\sup_{v \in f} \sP_v^\omega( H_f^+ < \bar H_U)$---we are in fact going to estimate the complementary probability $\inf_{v \in f} \sP_v^\omega( H_f^+ > \bar H_U)$.

Now intervenes the contraction of $f = \{x,y\}$ in the graph, following the procedure explained in the discussion of Theorem~\ref{contraction}. For the sake of simplicity, we will just estimate $\sP_x^\omega( H_f^+ > \bar H_U)$---the same reasoning, exchanging $x$ and~$y$, gives the same final estimate for $\sP_y^\omega( H_f^+ > \bar H_U)$. 

Consider the graph $\cG_U = (V_U, E_U)$ deduced from $\ZZ^d$ by replacing $\ZZ^d \setminus U$ by a cementery~$\partial$. Contract the edge $f = \{x, y\}$ into a single point $x_f$ and define a new environment over the new graph $\cG_U^f = (V_U^f, E_U^f)$ as follows: for every $\omega\in\Omega$, we define a corresponding environment $\tilde\omega$ over $\cG_U^f$ by setting, for every $(u, v) \in (V_U^f)^2$:
\begin{align*}
\tilde\omega(u, v) =
\begin{cases}
\omega(u, v) &\text{if } u \neq x_f \text{ and } v \neq x_f, \\
\omega(u, x) &\text{if } v = x_f \text{ and } (u, x) \in E_U, \\
\omega(u, y) &\text{if } v = x_f \text{ and } (u, y) \in E_U, \\
{s_f}{\omega(x, v)} &\text{if } u = x_f \text{ and } (x, v) \in E_U, \\
{s_f}{\omega(y, v)} &\text{if } u = x_f \text{ and } (y, v) \in E_U. \\
\end{cases}
\end{align*}

For every $\omega\in\Omega$, we have:
\begin{align*}
\sP_x^\omega( H_f^+ > \bar H_U)
&= \sum_{\sigma \in \Pi_x^f(U)} \omega_\sigma = \sum_{\sigma \in \Pi_x^f(U)} \omega_{\sigma_\mathrm{b/f}} \omega_{\sigma_\mathrm{odd}} \omega_{\sigma_\mathrm{exit}}
\end{align*}
where $\Pi_x^f(U) $ is the set of paths~$\sigma$ starting at~$v$ that exits~$f$ in finite time and then exits~$U$ without returning ever to~$f$, and where we have decomposed for every $\sigma\in\Pi_x^f(U)$ the product $\omega_\sigma = \prod_{e\in\sigma}\omega(e)$ in several factors:
\begin{itemize}
\item the part $\omega_{\sigma_\mathrm{b/f}}$ corresponding to the back-and-forths in $f$: if we sum the contributions of $\omega_{\sigma_\mathrm{b/f}}$ in $\omega_\sigma$, it introduces a factor $\sum_{i \pg 0} [\omega(x, y) \omega(y, x)]^i = \frac{1}{1 - p_f}$.
\item the weight of the eventual edge corresponding at the last crossing of the trap, if the time spent in~$f$ is odd : $\omega_{\sigma_\mathrm{odd}} = \omega(x, y) $ if this time is odd, $= 1$ otherwise. Note that, in both cases, because $f$ is a trap, $\omega_{\sigma_\mathrm{odd}} > \frac{1}{2}$.
\item the part $\omega_{\sigma_\mathrm{exit}}$ between the exit time from~$f$ and the exit from~$U$: summing this contribution, we get ${\sP_x^\omega(H_x^+ \wedge H_y^+ > \bar H_U)}$ or ${\sP_y^\omega(H_x^+ \wedge H_y^+ > \bar H_U)}$, depending on whether the ``odd" edge exists.
\end{itemize}
Therefore:
\begin{align*}
&\sP_x^\omega( H_f^+ > \bar H_U) \\
&\qquad = (1- p_f)^{-1} \left( \sP_x^\omega(H_x^+ \wedge H_y^+ > \bar H_U) + \omega(x,y) \sP_y^\omega(H_x^+ \wedge H_y^+ > \bar H_U) \right) \\
&\qquad\pg \frac{1}{1- p_f} \left( \frac{1}{2} \sP_x^\omega(H_x^+ \wedge H_y^+ > \bar H_U) +  \frac{1}{2} \sP_y^\omega(H_x^+ \wedge H_y^+ > \bar H_U) \right) \\
&\qquad\pg \frac{1}{2s_f(1- p_f)} \left( s_f \sP_x^\omega(H_x^+ \wedge H_y^+ > \bar H_U) + s_f \sP_y^\omega(H_x^+ \wedge H_y^+ > \bar H_U) \right) \\
&\qquad\pg \frac{1}{2} \left( s_f \sP_x^\omega(H_x^+ \wedge H_y^+ > \bar H_U) + s_f \sP_y^\omega(H_x^+ \wedge H_y^+ > \bar H_U) \right),
\end{align*}
where, for the last inequality, we used~\eqref{eq:pf-sf-general}.

We would like to replace $\sP^\omega$ by $\sP^{\tilde\omega}$. Note that, the difference between $\omega$ and~$\tilde\omega$, on the event that the walk~$X$ in the original graph (respectively in the contracted graph) starts from a vertex in~$f$ (respectively from~$x_f$) and never visits again a vertex in~$f$ (respectively~$x_f$), is only noticeable when crossing the first edge $(X_0, X_1)$, because, on this event, it is only for time~$k=1$ that $X_k \not\in f$ and $X_{k-1} \in f$ (respectively $X_k \neq x_f$ and $X_{k-1} = x_f$). This difference between $\omega$ and~$\tilde\omega$ is only the factor~$s_f$. We thus have:
$$ s_f \sP_x^\omega(H_x^+ \wedge H_y^+ > \bar H_U) + s_f \sP_y^\omega(H_x^+ \wedge H_y^+ > \bar H_U)
= \sP_{x_f}^{\tilde\omega}(H_{x_f}^+ > \bar H_U).$$
Finally, we just proved that:
$$ \sP_x^\omega( H_f^+ > \bar H_U) \pg \frac{1}{2} \sP_{x_f}^{\tilde\omega}(H_{x_f}^+ > \bar H_U). $$
The same bound is true for $\sP_x^\omega( H_f^+ > \bar H_U)$, and therefore:
$$ \sup_{v \in f} \sP_v^\omega( H_f^+ > \bar H_U) \pg \frac{1}{2} \sP_{x_f}^{\tilde\omega}(H_{x_f}^+ > \bar H_U). $$

Returning to our quantity of interest, we deduce that:
$$ \sE_0\left[(N_{f \vert U})^\beta\right]
\pp \EE^{(\alpha)}\left[\sum_{k=0}^{\infty} \left[(k+1)^\beta - k^\beta\right] \left[1 - \frac{1}{2} \sP_{x_f}^{\tilde\omega} ( H_{x_f}^+ > \bar H_U)\right]^k \right].
$$

Using another time Fubini's theorem, we recognize that this sum is in fact the $\beta$-th moment of a random variable which follows under~$\sP_{x_f}^{\tilde\omega}$ a geometric distribution with parameter $ \frac{1}{2} \sP_{x_f}^{\tilde\omega} ( H_{x_f}^+ > \bar H_U) $. Using lemma~\ref{appendix} (in Appendix), there exists a constant $C_\beta$ such that:
$$ \sE_0\left[(N_{f \vert U})^\beta\right]
\pp C_\beta \EE^{(\alpha)}\left[ \left(\frac{1}{2} \sP_{x_f}^{\tilde\omega} ( H_{x_f}^+ > \bar H_U)\right)^{-\beta} \right].
$$

The environment $\tilde\omega$ is not a Dirichlet environment, but lemma~\ref{contraction} ensures that, under $\sP_0^{\tilde\omega}$, there exists a constant $C_f$ such that:
$$ \EE^{(\alpha)} \left[ \sP_{x_f}^{\tilde\omega} ( H_{x_f}^+ > \bar H_U)^{-\beta} \right] \pp C_f \EE^{(\alpha_f)}\left[ \sP_{x_f}^{\omega} ( H_{x_f}^+ > \bar H_U)^{-\beta} \right], $$
where the weights $\alpha_f$ endowing the graph~$\cG_U^f$ are defined as follows:
\begin{align*}
\alpha_f(u, v) =
\begin{cases}
\alpha(u, v) &\text{if } u \neq x_f \text{ and } v \neq x_f, \\
\alpha(u, x) &\text{if } v = x_f \text{ and } (u, x) \in E_U, \\
\alpha(u, y) &\text{if } v = x_f \text{ and } (u, y) \in E_U, \\
\alpha(x, v) &\text{if } u = x_f \text{ and } (x, v) \in E_U, \\
\alpha(y, v) &\text{if } u = x_f \text{ and } (y, v) \in E_U. \\
\end{cases}
\end{align*}

In fact, reading carefully lemma~\ref{contraction}, we see that the constant $C_f$, which \textit{a priori} depends on the graph~$\cG_U^f$, and therefore on the choice of~$U$ and on the edge $f$, depends only on the number of contracted vertices, the degree in $\ZZ^d$ of the contracted vertices, the degree in $\cG_U$ of the new vertex after contraction and the weights of the edges issued from the contracted vertices. Note that the first three quantities do not vary with $U$ and~$f$, and that the last one only depends of the direction of~$f$. Therefore, there are only $d$ different $C_f$'s, and there exists a constant $C < \infty$ such that for every undirected edge~$f$ and every subset of vertices~$U$ containing~$f$:
$$ \sE_0\left[(N_{f \vert U})^\beta\right]
\pp 2^\beta C_\beta C \EE^{(\alpha_f)}\left[ \sP_{x_f}^{\omega}( H_{x_f}^+ > \bar H_U)^{-\beta} \right]. $$

In a similar fashion as what we have done in lemma~\ref{out_lemma1}, $\sP_{x_f}^{\omega}( H_{x_f}^+ > \bar H_U)$ can be stochastically lower-bounded by a Beta random variable, up to a change of measure from the Dirichlet environment with parameters~$\alpha_f$ to another Dirichlet environment. Precisely, we have, for every $\gamma>0$, $s>0$ and $p, q > 1$ conjugated exponents, and for every unitary flow $\theta_f$ from~$x_f$ to~$\partial$:
\begin{align*}
&\EE^{(\alpha_f)}\left[ \sP_{x_f}^{\omega} \left( H_{x_f}^+ > \bar H_U\right)^{-\beta} \right] \\
&\qquad\pp C_\beta (2d)^s \sum_{\underline{e} = x_f} \EE^{(\alpha_f)}[\text{Beta}(\gamma, \alpha_{x_f})^{-p\beta}]^{1/p} \frac{Z_{\alpha_f+\gamma\theta_f}}{Z_{\alpha_f}} \EE^{(\alpha_f+\gamma\theta_f)}[\omega^{-\gamma q \theta_f} \omega(e)^{qs}]^{1/q}.
\end{align*}

As in lemma~\ref{out_lemma2}, in the sum, the first expectation is bounded by a finite constant~$A$ independent of $f$ and $U$ as soon as $p\beta < \gamma$.

For the other factors in the sum, we can bound them by $A'$, with $A'$ an universal finite constant, as soon as we choose the parameters as in the recapitulation phase of lemma~\ref{out_lemma2}'s proof, with the sole difference that $\kappa'$ is replaced by $\min \kappa(S')$ over all finite set $S'$ constructed as follows: if $S$ is a connected subset of vertices of the contracted graph containing $x_f$ but not reduced to $x_f$, we denote $S'$ the set of vertices in the original graph deduces from $S$ by removing $x_f$ and adding instead $\overline{f}$ and $\underline{f}$. It is easy to prove that
$$ \sum_{e \in \partial_+ S} \alpha_f(e) = \sum_{e \in \partial_+ S'} \alpha(e) = \kappa(S'), $$
and, since $\{x, y\} \subsetneq S'$, $\kappa(S') > \kappa$. Again, $A$ and $a$ are independent on $f$ and $U$.

Therefore, we can let $U$ tend to $\ZZ^d$, and get by the monotonous convergence theorem:
$$ \sE_0^{(\alpha)}[(N_f)^\beta] = \lim_{U \to \ZZ^d} \sE_0^{(\alpha)}[(N_{f \vert U})^\beta] \pp C C_\beta (2d)^{s} (4d-2) A A'. $$

Theorem~\ref{Green_integr} concludes the proof.
\end{proof}

\subparagraph*{Remark 2.} At this point, we can see why our argument fails in dimension~2. For the same reason as in Remark~1, our argument in dimension ${d = 2}$ gives the following bound: there exists $\beta > \kappa$, $\beta' > 0$ and a constant $C_\beta < \infty$ such that:
$$ \sE_0\left[(N_f)^\beta \Ind{T_1 < \bar{H}_U} \right] < C_\beta \vert U \vert^{\beta'}. $$
This bound in $d = 2$ is not sufficient for the sequel of our argument. Compared to the case $d \pg 3$, it introduces a power of $\vert U \vert$, which is of order $\ln(x)$. Used in lemmas~\ref{traps_often_weak} and \ref{traps_often_strong}, it would introduce in their bounds the same logarithmic factor, which we cannot get rid of by adjusting $h(x)$ and~$m$ (see section~\ref{traps_often}).

\medskip

We can however improve lemma~\label{visits}, by controlling not only a single $N_f$, but the sum of all $N_f$'s for $f$ encountered before the first renewal time~$T_1$. In the following lemma, $\mathbf{(T)}$-condition is a crucial element.

\begin{lemma}
\label{visits_traps}
Assume $\mathbf{(T)}$-condition in $\ZZ^d, d \pg 3$. Then there exists $\beta \in (\kappa, 2)$ such that:
$$ \sE_0\left[\left(\sum_{f\in\cT_1} N_f \right)^\beta \right] < \infty. $$
\end{lemma}

\begin{proof}
Let $\beta \in (\kappa, 2)$. Since $\frac{\beta}{2} \in (0,1)$, we have:
\begin{align*}
\sE_0\left[\left(\sum_{f\in\cT_1} N_f \right)^\beta \right]
&= \sE_0\left[\left(\left(\sum_{f\in\cT_1} N_f \right)^{\beta/2}\right)^2 \right] \\
&\pp \sE_0\left[\left(\sum_{f\in\cT_1} N_f ^{\beta/2}\right)^2 \right] \\
&\pp \sum_{(f, f')\in\tilde{E}^2} \sE_0\left[N_f^{\beta/2} N_{f'}^{\beta/2}\right] \\
&\pp \sum_{(f, f')\in\tilde{E}^2} \sE_0\left[N_f^{\beta}\right]^{1/2} \sE_0\left[N_{f'}^{\beta}\right]^{1/2} \\*
&\pp \left(\sum_{f\in\tilde{E}} \sE_0\left[N_f^{\beta}\right]^{1/2}\right)^2.
\end{align*}

Fix an unoriented edge $f \in \tilde{E}$. Let $p>1$ a parameter to be chosen later and $q$ be its conjugate exponent. By Hölder's inequality:
\begin{align*}
\sE_0\left[N_f^{\beta}\right] &\pp \sE_0\left[(N_f)^{\beta} \Ind{N_f \pg 1} \right] \pp \sE_0\left[(N_f)^{p\beta}\right]^{1/p} \sP_0(N_f \pg 1)^{1/q}.
\end{align*}

By the previous lemma, we can choose $\beta > \kappa$ and $p$ such that $\sE_0\left[N_f^{p\beta}\right]$ is finite and bounded by a universal constant $C_{p\beta}$ for all unoriented edge~$f$. Therefore:
\begin{align*}
\sE_0\left[\left(\sum_{f\in\cT_1} N_f \right)^\beta \right]
&\pp C_{p\beta} \left(\sum_{f\in\tilde{E}} \sP_0(N_f \pg 1)^{1/2q} \right)^2.
\end{align*}

We thus need to show that this last sum is finite. Here is used $\mathbf{(T)}$-condition. Let $f$ be an unoriented edge, and denote by $\| f \|$ the norm of $f$ (that is to say, the minimal norm of its extremities). If $N_f \pg 1$, \textit{i.e.} if the random walk visits $f$ before $T_1$, necessarily the maximal distance $R_1$ reached by the walk before time $T_1$ is at least $\| f \|$. Therefore, by $\mathbf{(T)}$, there exists a positive constant $c$ such that, $f \in \tilde{E}$:
$$ \sP_0(N_f \pg 1) \pp \sP_0(R_1 \pg \| f \|) \pp e^{-c\| f \|}. $$

Finally:
\begin{align*}
\sum_{f\in\tilde{E}} \sP_0(N_f \pg 1)^{1/2q}
&\pp \sum_{n\pg 0} \sum_{\substack{f\in\tilde{E} \tq \\ \|f\|=n}} \sP_0(N_f \pg 1)^{1/2q} \\
&\pp \sum_{n\pg 0} \sum_{\substack{f\in\tilde{E} \tq \\ \|f\|=n}} e^{-\frac{c}{2q}n} \\
&\pp \sum_{n\pg 0} C n^{d-1} e^{-\frac{c}{2q}n} \\*
&< \infty,
\end{align*}
and:
$$\sE_0\left[\left(\sum_{f\in\cT_1} N_f \right)^\beta \right] < \infty.$$
\end{proof}

	\subsection{The time spent in the traps visited many times}

In general, the time spent in a trap can be dominated by a i.i.d. sum of geometric random variables, whose parameter is given by the strength of the trap and whose number is given by the number of visits to the trap. More precisely:
\begin{lemma}
Fix $f = \{x, y\}$ an unoriented edge. Then, knowing the number $N_f$ of visit to~$f$ before time~$T_1$ and the strength $s_f$ of~$f$, there exists a sequence of independent random variables $(H_f^j)_{j \pg 0}$, independent from $N_f$ and with geometric distribution of parameter $\frac{1}{2s_f}$, such that, $\sP_0$-almost surely:
$$ \ell_f^j \pp 2 H_f^j .$$
\end{lemma}

\begin{proof}
Note that, knowing the partially forgotten walk, the $\ell_f^j$, for $1 \pp j \pp N_f$, are independent random variables.

Knowing $\omega\in\Omega$, we can decompose $\ell_f^j$ as follows. First, we consider the time spent in performing back-and-forths in the trap. A back-and-forth occurs with probability $p_f = \omega(x,y) \omega(y,x)$. Therefore, the walk performs at least $k$ back-and-forths in~$f$ with probability $(p_f)^k$. Let $B_f^j$ the number of back-and-forths; from:
$$\PP(B_f^j \pg k) = (p_f)^k,$$
we deduce:
$$\PP(B_f^j = k-1) = (p_f)^{k-1} (1-p_f).$$
In other words, $B_f^j + 1$ follows a geometric distribution with parameter $1-p_f$. Secondly, there may remain to cross a last time the edge $f$, if $\ell_f^j$ is odd.

Therefore, $\ell_f^j \pp 2(B_f^j + 1)$, which means that $\ell_f^j$ is upper-bounded by twice a geometric distribution with parameter~$1-p_f$.

Remark that the geometric distributions with parameter $p > 0$ are stochastically decreasing when $p$ grows. We also know that, if $s_f \pg 2$, then:
$$ 1 - p_f \pg \frac{1}{2s_f}.$$
Therefore $\ell_f^j$ is stochastically dominated by twice a geometric distribution with parameter $\frac{1}{2s_f}$.
\end{proof}

We now address separately the cases whether the strength of the trap is lower-bounded or upper-bounded.

\subsubsection{Traps whose strength is upper-bounded}

The previous lemma enables us to control the time spent in a given trap, if we know that its strength is at most $h$.

\begin{lemma}
Let $f$ be an unoriented edge, $\beta > \kappa$, $\eta > 0$, and $h \pg 0$. Then there exists a finite constant $C_{\beta, \eta}$ such that:
$$ \sE_0\left[N_f^\eta \left(\sum_{j=1}^{N_f} \ell_f^j\right)^{\beta} \Ind{2 \pp s_f \pp h}\right] \pp C_{\beta, \eta} \sE_0[{N_f}^{\beta+\eta}] h^{\beta-\kappa}.$$
\end{lemma}

\begin{proof}
Fix $f\in\tilde{E}$, $\beta > \kappa$, $\eta > 0$, and $h \pg 0$. We will first consider:
\begin{align*}
&\sE_0\left[\left. N_f^\eta \left(\sum_{j=1}^{N_f} \ell_f^j\right)^{\beta} \Ind{2 \pp s_f \pp h} \right\vert N_f, s_f\right] \\
&\qquad\pp \sE_0\left[\left.N_f^{\eta+\beta-1} \sum_{j=1}^{N_f} \left(\ell_f^j\right)^{\beta} \right\vert N_f, s_f\right] \Ind{2 \pp s_f \pp h} \\
&\qquad\pp N_f^{\eta+\beta-1} \sum_{j=1}^{N_f} \sE_0\left[\left.\left(\ell_f^j\right)^{\beta} \right\vert N_f, s_f\right] \Ind{2 \pp s_f \pp h} \\
&\qquad\pp N_f^{\eta+\beta-1} \sum_{j=1}^{N_f} \sE_0\left[\left.\left(H_f^j\right)^{\beta} \right\vert N_f, s_f\right] \Ind{2 \pp s_f \pp h)} \\
&\qquad\pp C_\beta N_f^{\eta+\beta} s_f^{\beta} \Ind{2 \pp s_f \pp h}
\end{align*}
where the $H_f^j$ are the random variables introduced in the previous lemma, and where the last inequality comes from lemma~\ref{appendix}. Integrating over $N_f$ and $s_f$, we get:
\begin{align*}
\sE_0\left[N_f^\eta \left(\sum_{j=1}^{N_f} \ell_f^j\right)^{\beta} \Ind{2 \pp s_f \pp h}\right]
\pp C_\beta \sE_0\left[N_f^{\beta+\eta} s_f^{\beta} \Ind{2 \pp s_f \pp h}\right].
\end{align*}

We now want to use the quasi-independence between $N_f$ and $s_f$ (lemma~\ref{quasi-independence}). It can be done by replacing in the previous inequality $s_f^\beta$ by $\int_0^{s_f^\beta} \dd t$. Indeed:
\begin{align*}
&\sE_0\left[N_f^\eta \left(\sum_{j=1}^{N_f} \ell_f^j\right)^{\beta} \Ind{2 \pp s_f \pp h}\right] \\
&\qquad\pp C_\beta \sE_0\left[N_f^{\beta+\eta} \left(\int_0^{s_f^{\beta}} \dd t\right) \Ind{2 \pp s_f \pp h}\right] \\
&\qquad\pp C_\beta \sE_0\left[N_f^{\beta+\eta} \int_0^{\infty} \Ind{t^{1/\beta} \pp s_f \pp h, s_f \pg 2} \dd t \right] \\
&\qquad\pp C_\beta \sE_0\left[N_f^{\beta+\eta}  \left(2^{\beta} + \int_{2^\beta}^{h^\beta} \Ind{s_f \pg t^{1/\beta}} \dd t \right)\right] \\
&\qquad\pp C_\beta \left( 2^{\beta} \sE_0\left[N_f^{\beta+\eta}\right] + \int_{2^\beta}^{h^\beta} \sE_0\left[N_f^{\beta+\eta} \Ind{s_f \pg t^{1/\beta}}\right] \dd t \right).
\end{align*}

By lemma~\ref{quasi-independence}, if we choose $\beta$ and $\eta$ such that $\beta \in [1,2)$ and $\beta + \eta \in (1,2)$, so that $x \longmapsto x^{\beta+\eta}$ has a concave increasing derivative, there exists a finite constant $C_{\beta+\eta}$ such that
\begin{align*}
\sE_0\left[N_f^\eta \left(\sum_{j=1}^{N_f} \ell_f^j\right)^{\beta} \Ind{2 \pp s_f \pp h}\right]
&\pp C_\beta \sE_0\left[N_f^{\beta+\eta}\right] \left( 2^{\beta} + C_{\beta+\eta} \int_{2^\beta}^{h^\beta} t^{-\kappa/\beta} \dd t \right) \\
&\pp C_{\beta, \eta} \sE_0[{N_f}^{\beta+\eta}] h^{\beta-\kappa}.
\end{align*}
\end{proof}

The previous lemma enables us to control the time spent in all traps whose strength is not too big and which are visited many times:

\begin{lemma}
\label{traps_often_weak}
Let $m \in \NN^*$, and $h \pg 0$. If $\kappa \in (1,2)$ and $\mathbf{(T)}$-condition holds, then there exists $\beta \in (\kappa, 2)$, $\eta > 0$ and a finite constant $C_{\beta, \eta}$ such that:
$$ \sP_0\left(\sum_{f \in \cT_1} \sum_{j=1}^{N_f} \ell_f^j \Ind{N_f \pg m, 2 \pp s_f \pp h} \pg x \right) \pp C_{\beta, \eta} m^{-\eta} x^{-\beta} h^{\beta-\kappa}. $$
\end{lemma}

\begin{proof}
Let $m \in \NN^*$, and $h \pg 0$. Let $\beta \in (\kappa, \beta)$ and $\eta > 0$. By Markov's inequality:
\begin{align*}
&\sP_0\left(\sum_{f \in \cT_1} \sum_{j=1}^{N_f} \ell_f^j \Ind{N_f \pg m, 2 \pp s_f \pp h} \pg x \right) \\
&\qquad \pp x^{-\beta} \sE_0\left[\left(\sum_{f \in \tilde \cT_1} \sum_{j=1}^{N_f} \ell_f^j \Ind{N_f \pg m, 2 \pp s_f \pp h} \right)^\beta \right].
\end{align*}

We then proceed as in lemma~\ref{visits_traps} and use another time Markov's inequality:
\begin{align*}
&\sP_0\left(\sum_{f \in \cT_1} \sum_{j=1}^{N_f} \ell_f^j \Ind{N_f \pg m, 2 \pp s_f \pp h} \pg x \right) \\
&\qquad \pp x^{-\beta} \left(\sum_{f \in \tilde E} \sE_0\left[ \left(\sum_{j=1}^{N_f} \ell_f^j\right)^\beta \Ind{N_f \pg m, 2 \pp s_f \pp h} \right]^{1/2} \right)^2 \\
&\qquad \pp x^{-\beta} m^{-\eta} \left(\sum_{f \in \tilde E} \sE_0\left[ N_f^{\eta} \left(\sum_{j=1}^{N_f} \ell_f^j\right)^\beta \Ind{2 \pp s_f \pp h} \right]^{1/2} \right)^2.
\end{align*}

By the previous lemma, we get:
\begin{align*}
&\sP_0\left(\sum_{f \in \cT_1} \sum_{j=1}^{N_f} \ell_f^j \Ind{N_f \pg m, 2 \pp s_f \pp h} \pg x \right) \\
&\qquad\pp C_{\beta, \eta} x^{-\beta} m^{-\eta} h^{\beta-\kappa} \times \left(\sum_{f \in \tilde E} \sE_0[{N_f}^{\beta+\eta}]^{1/2} \right)^2.
\end{align*}
If we choose $\beta$ and $\eta$ such that $\sE_0[{N_f}^{\beta+\eta}]<\infty$ for every $f$, lemma~\ref{quasi-independence} ensures us that the last sum is finite, which proves the announced result.
\end{proof}

\subsubsection{Traps whose strength is lower-bounded}
	
	We now turn to the probability that the walk visits many time a trap of strength at least $h$ before time $T_1$. Precisely:
	
\begin{lemma}
\label{traps_often_strong}
Under $\mathbf{(T)}$-condition, there exists $\eta > \kappa$ and a finite constant $C_\eta$ such that, for every $m \in \NN^*$ and $h$ large enough,
$$\sP_0\left(\exists f \in \cT_1 \tq N_f \pg m, s_f \pg h \right) \pp C_\eta m^{-\eta} h^{-\kappa}. $$
\end{lemma}

\begin{proof}
For every $\eta > 1$, by Markov's inequality:
\begin{align*}
\sP_0\left(\exists f \in \cT_1 \tq N_f \pg m, s_f \pg h \right)
&\pp \sP_0\left(\sum_{f\in\cT_1} N_f^\eta \Ind{s_f \pg h} \pg m^{\eta} \right) \\*
&\pp m^{-\eta} \sE_0\left[\sum_{f \in \cT_1} N_f^\eta \Ind{s_f \pg h} \right]
\end{align*}
and by lemma~\ref{quasi-independence}
\begin{align*}
\sP_0\left(\exists f \in \cT_1 \tq N_f \pg m, s_f \pg h \right)
&\pp m^{-\eta} \sE_0\left[\sum_{f \in \cT_1} N_f^\eta \right] h^{-\kappa} \\
&\pp C_\eta m^{-\eta} h^{-\kappa},
\end{align*}
the constant $C_\eta$ being given by lemma~\ref{visits_traps}.
\end{proof}

\subsubsection{Total time in traps visited many times}
\label{traps_often}

Let us sum up this section. Fix $\epsilon > 0$, and set, for $x > 0$:
\begin{align*}
h(x) &= \epsilon x.
\end{align*}

We can distinguish, within the whole time spent before $T_1$ in traps visited more than $m$~times, the two following durations:
\begin{enumerate}
\item The time spent in traps whose strength is at most $h(x)$. According to lemma~\ref{traps_often_weak}, applied for $m=1$ and $\eta=0$, there exists $\beta \in (\kappa, 2)$ such that the tail of this duration is at most $C_{\beta, \eta} \epsilon^{\beta-\kappa} x^{-\kappa}$. Note that $\beta > \kappa$, so that $\epsilon^{\beta-\kappa} \xrightarrow[\epsilon\to 0]{} 0$. 
\item The time spent in traps whose strength is at least $h(x)$. There exists some $\eta \in (\kappa, 2)$, such that, for every $m \in \NN^*$, the probability that such a trap exists is at most $C_\eta m^{-\eta} \epsilon^{-\kappa} x^{-\kappa}$, according to lemma~\ref{traps_often_strong}. 
Here we adjust $m \in \NN^*$ in terms of $\epsilon$ later, as announced in the sketch of proof (section~\ref{sketch}). Choose now $m = m(\epsilon) = \left\lfloor \epsilon^{\frac{-\kappa- 1}{\eta}} \right\rfloor$, so that the previous bound is at most of order $\epsilon x^{-\kappa}$. Therefore, the tail of the time spent in traps whose strength is at least $h(x)$ is at most of order $\epsilon x^{-\kappa}$.

\end{enumerate}
We easily deduce of these two tail estimates that for every $\epsilon > 0$, for $x$ large enough:
$$\sP_0\left(\sum_{f\in\cT_1} \sum_{j=1}^{N_f} \ell_f^j \Ind{N_f \pg m(\epsilon)} \pg x\right) \pp \frac{1}{2}(\epsilon + \epsilon^{\beta-\kappa}) x^{-\kappa}.$$

	\subsection{The time spent in traps visited a few times}

	We know turn to the estimate of:
	$$\sP_0\left(\sum_{f\in\cT_1} \sum_{j=1}^{N_f} \ell_f^j \Ind{N_f \pp m(\epsilon)} \pg x\right).$$

\subsubsection{Weak traps}

	By lemma~\ref{traps_often_weak}, in the special case $m = 1$ and $h = h(x) = \epsilon x$, we know that there exists $\beta > \kappa$ and $C_\beta < \infty$ such that, for every $\epsilon > 0$ and for $x$ large enough:
	$$\sP_0\left(\sum_{f\in\cT_1} \sum_{j=1}^{N_f} \ell_f^j \Ind{2 \pp s_f \pp \epsilon x} \pg x\right) \pp C_\beta \epsilon^{\beta-\kappa} x^{-\kappa}.$$
	\emph{A fortiori}, we have (here, as in the previous section, $m(\epsilon) = \left\lfloor \epsilon^{\frac{-\kappa- 1}{\eta}} \right\rfloor$):
	
\begin{lemma}
$\displaystyle \sP_0\left(\sum_{f\in\cT_1} \sum_{j=1}^{N_f} \ell_f^j \Ind{2 \pp s_f \pp \epsilon x, N_f \pp m(\epsilon)} \pg x\right) \pp C_\beta \epsilon^{\beta-\kappa} x^{-\kappa}.$
\end{lemma}
	
	We can thus focus on traps whose strength is lower-bounded by $h(x) = \epsilon x$.

\subsubsection{There is a single trap that significantly delays the walk}

	In fact, with high probability, there is only one trap whose strength is lower-bounded by $h(x) = \epsilon x$.
	
\begin{lemma}
	There exists a constant $C < \infty$ such that, for every $m > 0$ and $x$ large enough,
$$
\sP_0\left(\exists f, f' \in \cT_U^2 \tq N_f, N_{f'} \pp m; s_f, s_{f'} \pg h(x) \right)
\pp C |U|^2  \exp\left(\frac{5 m}{h(x)}\right) m^{2\kappa} h(x)^{-2\kappa}. 
$$
\end{lemma}

\begin{proof}
Remind that, knowing the partially forgotten environment and walk (and especially knowing the $N_f$'s), the strengths of the different traps are independent (lemma~\ref{density_forgotten}). Therefore, for every unoriented edges $f$ and $f'$:
\begin{align*}
&\sP_0\left(\left.N_f \pp m, N_{f'} \pp m, s_f \pg h(x), s_{f'} \pg h(x) \right| \tilde\omega, \tilde X \right) \\*
&\qquad = \sP_0\left(\left. s_f \pg h(x), s_{f'} \pg h(x) \right| \tilde\omega, \tilde X \right) \Ind{N_f \pp m, N_{f'} \pp m} \\
&\qquad\pp\sP_0\left(\left. s_f \pg h(x) \right| \tilde\omega, \tilde X \right) \Ind{N_f \pp m} \cdot \sP_0\left(\left. s_{f'} \pg h(x) \right| \tilde\omega, \tilde X \right) \Ind{N_{f'} \pp m}.
\end{align*}

Thanks to lemma~\ref{quasi-independence}, we have:
$$
\PP\left(\left. s_f \pg h(x) \right\vert \tilde\omega, \tilde{X} \right)
\pp C \exp\left(\frac{5 N_f}{2 h(x)}\right) (N_f)^{\kappa_j} h(x)^{-\kappa},
$$
so that:
\begin{align*}
&\sP_0\left(N_f \pp m, N_{f'} \pp m, s_f \pg h(x), s_{f'} \pg h(x) | \tilde\omega, \tilde X \right) \\
&\qquad \pp C \exp\left(\frac{5 (N_f+N_{f'})}{2h(x)}\right) (N_f N_f')^{\kappa_j} h(x)^{-2\kappa} \Ind{N_f \pp m, N_{f'} \pp m},
\end{align*}
and:
\begin{align*}
\sP_0\left(N_f \pp m, N_{f'} \pp m, s_f \pg h(x), s_{f'} \pg h(x) \right) &\pp C \exp\left(\frac{5 m}{h(x)}\right) m^{2\kappa_j} h(x)^{-2\kappa}.
\end{align*}

A union bound finally gives the announced estimate.
\end{proof}

\subsubsection{Strong traps}
	
	Therefore, with very high probability, the sum $\sum_{f \in \cT_1}\sum_{j=1}^{N_f} \ell_f^j \Ind{s_f > h(x), N_f \pp m(\epsilon)}$ is reduced to a single term. For this term, we can use the equivalent established in lemma~\ref{tail}, which gives exactly the order $x^{-\kappa}$ we expect:
	
\begin{lemma}
\label{equivalent}
	Let $\epsilon > 0$. With $h(x) = \epsilon x$ and $m(\epsilon) = \left\lfloor \epsilon^{\frac{-\kappa- 1}{\eta}} \right\rfloor$ as above, there exists a constant $C(\epsilon) < \infty$ such that, for every $x$ large enough,
$$
\sP_0\left(\sum_{f \in \cT_1} \sum_{j=1}^{N_f} \ell_f^j \Ind{s_f > h(x), N_f \pp m(\epsilon)} \pg x, T_1 < \bar{H}_U \right) \sim C(\epsilon) x^{-\kappa}. 
$$
Moreover, $C(\epsilon)$ converges to a finite constant $C_0$ when $\epsilon \to 0$.
\end{lemma}

\begin{proof}
For $m > 0$, denote by $\mathfrak{C}_m$ the set of possibles configurations $\mathfrak{c} = (j, N_1, N_2, N_3, N_4)$ such that $N_1 + \ldots + N_4 \pp m$. (Remind the definition of a configuration of a trap in section~\ref{traps}.)

By the previous lemma,
\begin{align*}
&\sP_0\left(\exists f, f' \in \cT_U^2 \tq N_f, N_{f'} \pp m(\epsilon); s_f, s_{f'} \pg h(x) \right) \\
&\qquad\pp C |U|^2  \exp\left(\frac{5 m(\epsilon)}{\epsilon x}\right) m(\epsilon)^{2\kappa} \epsilon^{-2\kappa} x^{-2\kappa} \\
&\qquad\pp C(\epsilon) |U|^2 x^{-2\kappa}
\end{align*}
Therefore, for all $\epsilon > 0$,
$$ x^{\kappa} \sP_0\left(\exists f, f' \in \cT_U^2 \tq N_f, N_{f'} \pp m(\epsilon); s_f, s_{f'} \pg h(x) \right) \xrightarrow[x\to+\infty]{} 0. $$

Let us so concentrate on the case where there is a single trap of strength at least~$h(x)$ and visited at most $m$~times. We clearly have:
\begin{align*}
&\sP_0\left(\exists ! f \in \cT_1 \tq s_f > h(x), 1 \pp N_f \pp m(\epsilon) \text{ and } \sum_{j=1}^{N_f} \ell_f^j \Ind{s_f > h(x), N_f \pp m(\epsilon)} \pg x, T_1 < \bar{H}_U \right) \\*
&= \sum_{f \in \tilde E_U} \sP_0\left(\left\lbrace f \text{ is the unique trap such that } s_f > h(x), 1 \pp N_f \pp m(\epsilon)\right\rbrace \cap \left\lbrace\sum_{j=1}^{N_f} \ell_f^j \pg x \right\rbrace\right).
\end{align*}

Therefore, we have the following bounds:
\begin{align*}
&\sP_0\left(\exists ! f \in \cT_1 \tq s_f > h(x), 1 \pp N_f \pp m(\epsilon) \text{ and } \sum_{j=1}^{N_f} \ell_f^j \Ind{s_f > h(x), N_f \pp m(\epsilon)} \pg x, T_1 < \bar{H}_U \right) \\
&\pp \sum_{f \in \tilde E} \sP_0\left(\left\lbrace s_f > h(x), 1 \pp N_f \pp m(\epsilon)\right\rbrace \cap \left\lbrace\sum_{j=1}^{N_f} \ell_f^j \pg x \right\rbrace\right)
\end{align*}
and
\begin{align*}
&\sP_0\left(\exists ! f \in \cT_1 \tq s_f > h(x), 1 \pp N_f \pp m(\epsilon) \text{ and } \sum_{j=1}^{N_f} \ell_f^j \Ind{s_f > h(x), N_f \pp m(\epsilon)} \pg x, T_1 < \bar{H}_U \right) \\
&\pg \sum_{f \in \tilde E} \sP_0\left(\left\lbrace s_f > h(x), 1 \pp N_f \pp m(\epsilon)\right\rbrace \cap \left\lbrace\sum_{j=1}^{N_f} \ell_f^j \pg x \right\rbrace\right) \\
&\qquad - \sum_{f \in \tilde E} \sum_{f' \in \tilde E} \sP_0\left(\left\lbrace s_f, s_{f'} > h(x); 1 \pp N_f, N_{f'} \pp m(\epsilon)\right\rbrace \cap \left\lbrace\sum_{j=1}^{N_f} \ell_f^j \pg x \right\rbrace\right)
\end{align*}
Note that, in the second sum, the probabilities are at most $C(\epsilon) x^{-2k}$, and that the sum contains a number of terms of order $|U|^2$. Remind that $U$ is a box whose edge has length $\frac{1+\kappa}{c_{\mathrm{T}}} \ln(x)$. Therefore, for every $\epsilon > 0$, this second sum, multiplied by $x^{\kappa}$, vanishes when $x \to \infty$.

Let us so concentrate on:
\begin{align*}
&\sum_{f \in \tilde E} \sP_0\left( s_f > h(x), 1 \pp N_f \pp m(\epsilon), \sum_{j=1}^{N_f} \ell_f^j \pg x\right) \\
&\qquad = \sum_{\mathfrak{c} \in \mathfrak{C}_{m(\epsilon)}} \sum_{f \in \tilde E} \sP_0\left(\mathfrak{c}(f) = \mathfrak{c}, s_f > h(x), 1 \pp N_f \pp m(\epsilon), \sum_{j=1}^{N_f} \ell_f^j \pg x\right) \\
&\qquad = \sum_{\mathfrak{c} \in \mathfrak{C}_{m(\epsilon)}} \sum_{f \in \tilde E} \sE_0\left[ \sP_0\left(\left.\mathfrak{c}(f) = \mathfrak{c}, s_f > h(x), 1 \pp N_f \pp m(\epsilon), \sum_{j=1}^{N_f} \ell_f^j \pg x \right\vert \tilde X, \tilde \omega\right)\right] \\
&\qquad = \sum_{\mathfrak{c} \in \mathfrak{C}_{m(\epsilon)}} \sum_{f \in \tilde E} \sE_0\left[ \Ind{\mathfrak{c}(f) = \mathfrak{c}, 1 \pp N_f \pp m(\epsilon)} \sP_0\left(\left.s_f > h(x), \sum_{j=1}^{N_f} \ell_f^j \pg x \right\vert \tilde X, \tilde \omega\right)\right] \\
&\qquad \sim \sum_{\mathfrak{c} \in \mathfrak{C}_{m(\epsilon)}} \sum_{f \in \tilde E} \sE_0\left[ \Ind{\mathfrak{c}(f) = \mathfrak{c}, 1 \pp N_f \pp m(\epsilon)} C_{\mathfrak{c}} x^{-\kappa} \right] \\
&\qquad \sim \left( \sum_{\mathfrak{c} \in \mathfrak{C}_{m(\epsilon)}} C_{\mathfrak{c}}(\epsilon) \sum_{f \in \tilde E} \sP_0\left( \mathfrak{c}(f) = \mathfrak{c}, 1 \pp N_f \pp m(\epsilon) \right) \right) x^{-\kappa}
\end{align*}
where we used in the penultimate line lemma~\ref{tail}. It suffices to call $C(\epsilon)$ this last sum---it is finite for every~$\epsilon > 0$ since $\mathfrak{C}_{m(\epsilon)}$ is a finite set.

It remains to note that $C(\epsilon)$ converges when $\epsilon$ goes to~0. Because every $C_{\mathfrak{c}}(\epsilon)$ grows as $\epsilon$ goes to~0, so does $C(\epsilon)$: it therefore tends to some (possibly infinite) constant $C_0$; but, by the same reasoning, $C_0 - C(\epsilon)$ controls in fact the tail of the time spent in the traps visited more than $m(\epsilon)$ times, which we proved to be negligible with respect to $x^{-\kappa}$ in section~\ref{traps_often} this compels $C_0$ to be finite.
\end{proof}

	\subsection{Conclusion}
	
\begin{lemma}
\label{concl_lemma}
Under the assumptions of Theorem~\ref{main}, there exists a constant $C > 0$ such that:
$$ \sP_0(T_1 \pg x) \sim C x^{-\kappa}.$$
\end{lemma}

\begin{proof}
Let $x > 0$ be large enough.

\paragraph{Step 1.} We can assume furthermore that, between $0$ and $T_1$, the walk is contained in a box~$U(x)$ with edge of order $\frac{1+\kappa}{c_\mathrm{T}}\ln(x)$ (where $c_\mathrm{T}$ is the constant in $\mathbf{(T)}$).

Indeed, we clearly have:
$$ \sP_0(T_1 \pg x, T_1 < \bar{H}_{U(x)}) \pp \sP_0(T_1 \pg x) \pp \sP_0(T_1 \pg x, T_1 < \bar{H}_{U(x)}) + \sP_0(T_1 > \bar{H}_{U(x)}).$$
By $\mathbf{(T)}$:
\begin{align*}
\sP_0(T_1 > \bar{H}_{U(x)})
&\pp \sP_0\left( \max_{1 \pp n \pp T_1} \|X_n\| \pg \frac{1+\kappa}{c_\mathrm{T}}\ln(x) \right) \\
&\pp \sE_0\left[ e^{c_\mathrm{T} \max_{1 \pp n \pp T_1} \|X_n\|} \right] e^{-(\kappa+1) \ln(x)} \\
&= \grando{x^{-\kappa -1}} = \petito{x^{-\kappa}}.
\end{align*} 

Therefore, it suffices to show that:
$$ \sP_0(T_1 \pg x, T_1 < \bar{H}_{U(x)}) \sim C x^{-\kappa}.$$

\paragraph{Step 2.} We can neglect the time spent out of traps. In other words, using notations of section~\ref{traps}, it suffices to show that:
$$ \sP_0(T_1^\bullet \pg x, T_1 < \bar{H}_{U(x)}) \sim C x^{-\kappa} $$
to be ensured that:
$$ \sP_0(T_1 \pg x, T_1 < \bar{H}_{U(x)}) \sim C x^{-\kappa}. $$

Indeed, assume that:
$$ \sP_0(T_1^\bullet \pg x, T_1 < \bar{H}_{U(x)}) \sim C x^{-\kappa}. $$

On the one hand, we clearly have:
$$\sP_0(T_1 \pg x, T_1 < \bar{H}_{U(x)}) \pg \sP_0(T_1^\bullet \pg x, T_1 < \bar{H}_{U(x)})$$
so:
$$ \liminf_{x \to \infty} x^\kappa \sP_0(T_1 \pg x, T_1 < \bar{H}_{U(x)}) \pg C. $$

On the other hand, let $\epsilon > 0$. We know by lemma~\ref{out_lemma2}, that there exists $\eta > 0$ such that:
$$ \sP_0(T_1^\circ \pg \delta x, T_1 < \bar{H}_{U(x)}) = \grando{x^{-\kappa-\eta}}. $$
Let $x_0 > 0$ such that for $x > \frac{x_0}{2}$:
$$ \sP_0(T_1^\bullet \pg x, T_1 < \bar{H}_{U(x)}) \pp (1+\epsilon) C x^{-\kappa}. $$

For every $x \pg x_0$ and $\delta > 0$ small enough such that $\delta < \frac{\epsilon}{8}$ and $(1-\delta)^{-\kappa} \pp 1 + 2\kappa\delta \pp 1 + \frac{\epsilon}{2} $  and:
\begin{align*}
\sP_0(T_1 \pg x, T_1 < \bar{H}_{U(x)})
&\pp \sP_0(T_1^\bullet \pg (1-\delta)x, T_1 < \bar{H}_{U(x)}) + \sP_0(T_1^\circ \pg \delta x, T_1 < \bar{H}_{U(x)}) \\
&\pp (1+\epsilon) (1-\delta)^{-\kappa} C x^{-\kappa} + \delta^{-\kappa-\eta} x^{-\kappa-\eta} \\
&\pp (1+\epsilon) (1+2\kappa\delta) x^{-\kappa} + \delta^{-\kappa-\eta} x_0^{-\eta} x^{-\kappa} \\
&\pp (1+\epsilon) \left(1+2\kappa\delta+\frac{1}{C \delta^{\kappa+\eta} x_0^{\eta}}\right) C x^{-\kappa} \\
&\pp (1+\epsilon)^2 C x^{-\kappa}
\end{align*}
provided that $x_0^{-\eta} \pp \frac{C}{2} \delta^{\kappa+\eta}$. This implies that
$$ \limsup_{x \to \infty} x^\kappa \sP_0(T_1 \pg x, T_1 < \bar{H}_{U(x)}) \pp (1+\epsilon)^2 C, $$
so letting $\epsilon \to 0$ gives:
$$ \limsup_{x \to \infty} x^\kappa \sP_0(T_1 \pg x, T_1 < \bar{H}_{U(x)}) \pp C. $$

The two limits coincide to show the desired equivalent.

It thus suffices to show that:
$$ \sP_0(T_1^\bullet \pg x, T_1 < \bar{H}_{U(x)}) \sim C x^{-\kappa}. $$

\paragraph{Step 3.} The same way as for the time spent out of traps, we can show that we can discard the time spent in weak traps or in strong traps visited many times, because, thanks to section~\ref{traps_often} and lemma~21, we know that there exists some $\beta \in (\kappa, 2)$ such that:
$$ \sP_0\left(\sum_{f\in\cT_1} \sum_{j=1}^{N_f} \ell_f^j \Ind{N_f \pg m(\epsilon)} \pg x\right) \pp \frac{1}{2}(\epsilon + \epsilon^{\beta-\kappa}) x^{-\kappa},$$
$$ \displaystyle \sP_0\left(\sum_{f\in\cT_1} \sum_{j=1}^{N_f} \ell_f^j \Ind{2 \pp s_f \pp \epsilon x, N_f \pp m(\epsilon)} \pg x, T_1 < \bar{H}_{U(x)} \right) \pp C_\beta \epsilon^{\beta-\kappa} x^{-\kappa}.$$

It thus suffices to prove that, for every $\epsilon > 0$, that:
\begin{align*}
C^-(\epsilon) &\pp \liminf_{x\to\infty} x^\kappa \sP_0\left(\sum_{f\in\cT_1} \sum_{j=1}^{N_f} \ell_f^j \Ind{s_f \pg h(x), N_f \pp m(\epsilon)} \pg x, T_1 < \bar{H}_{U(x)}\right) \\*
&\pp \limsup_{x\to\infty} x^\kappa \sP_0\left(\sum_{f\in\cT_1} \sum_{j=1}^{N_f} \ell_f^j \Ind{s_f \pg h(x), N_f \pp m(\epsilon)} \pg x, T_1 < \bar{H}_{U(x)}\right) \pp C^+(\epsilon)
\end{align*}
where $C^\pm(\epsilon)$ converge to a constant $C$ when $\epsilon \to 0$, to be ensured that:
$$ \sP_0(T_1^\bullet \pg x, T_1 < \bar{H}_{U(x)}) \sim C x^{-\kappa}. $$

But these inequalities are a direct consequence of lemma~\ref{equivalent}. This concludes the proof.
\end{proof}
	
\section{Limit theorem for the random walk}

We now turn to the proof of our main result, Theorem~\ref{main}.

\begin{proof}[Proof of Theorem~\ref{main}] The proof consists in applying the stable law theorem to the renewal times so as to get the distribution of their fluctuations, and then to apply common inversion techniques so as to get the distributions of the fluctuations of the random walk itself. We mainly adopt the notations of~\cite{bowditch2021fluctuations}.

\paragraph{Step 1: limit theorem for $T$.}
We know that the renewal intervals $(T_{i+1} - T_i)_{i \pg 1}$ are i.i.d. random variables. By lemma~\ref{concl_lemma}, their tail is of order $x^{-\kappa}$ exactly. Therefore, by Theorem 4.5.3 of~\cite{whitt2002limits}, we have, in the $J_1$-topology:
\[
\left(t \longmapsto n^{-\frac{1}{\kappa}}\left(T_{\lfloor nt \rfloor} - n \tau t \right)\right)\xrightarrow[n\to\infty]{\sP_0^{(\alpha)}\text{-(d)}} c_1 \mathcal{S}^{\kappa}
\]
where we have denoted $\tau = \sE_0^{(\alpha)}[T_2 - T_1]$.

\paragraph{Step 2: limit theorem for $\Delta$.} Set, for $a \pg 0$ and $k \in \NN$,
$$\Delta_{a} = \inf\ens{i \pg 0 \tq \langle X_i, \hat u \rangle \pg a}.$$

Define, for $k \in \NN$, the $k$-th renewal level:
$$\rho_k = \langle X_{T_k}, \hat u \rangle$$
and denote by $\psi$ its càdlàg inverse: for all $a > 0$,
$$\psi(a) = \inf\ens{k \pg 0 \tq \rho_k \pg a}.$$

By definition of $\psi$, for every $a > 0$,
$$ T_{\psi(na) - 1} \pp \Delta_{na} \pp T_{\psi(na)}.$$
Therefore, before considering the limit of the process $\left(a \longmapsto n^{-\frac{1}{\kappa}}\left(\Delta_{na} - n \langle v, \hat u \rangle^{-1} a \right)\right)$, we determine the limit of the process $\left(a \longmapsto n^{-\frac{1}{\kappa}}\left(T_{\psi(na)} - n \langle v, \hat u \rangle^{-1} a \right)\right)$. Decompose for every $a > 0$:
\begin{align*}
T_{\psi(na)} - n \langle v, \hat u \rangle^{-1} a
= \phantom{+\:} & [T_{\psi(na)} - \tau\psi(na)] \\*
+\: &{}[\psi(na)\tau - n \langle v, \hat u \rangle^{-1} a].
\end{align*}
The first square bracket reminds us of Step~1, and the second one prompts us into finding the fluctuations of $\psi$. Let us so determine the asymptotic behaviour of $\psi$. As $\psi$ is the càdlàg inverse of $\rho$, this can be deduced from the asymptotic behaviour of $\rho$.

Note that, by condition $\mathbf{(T)}$, the $\rho_{k+1} - \rho_k$ ($k \pg 1$), which are i.i.d. random variables with expectation $\tau \langle v, \hat u \rangle$ (where $v$ is the asymptotic velocity of the random walk), have exponential moments. As a consequence, they satisfy a law of large numbers: in $J_1$-topology,
\[
\left(t \longmapsto n^{-1} \rho_{\lfloor nt \rfloor}\right) \xrightarrow[n\to\infty]{} \left(t \longmapsto \tau \langle v, \hat u \rangle t \right).
\]
Since this limit is a strictly increasing function, by Corollary 13.6.4 of \cite{whitt2002limits}: in $J_1$-topology,
\[
\left(a \longmapsto n^{-1} \psi(na)\right) \xrightarrow[n\to\infty]{} \left(a \longmapsto (\tau \langle v, \hat u \rangle)^{-1} a \right),
\]
and by continuity of composition (Theorem~13.2.2 of \cite{whitt2002limits}): in $J_1$-topology,
\[
\left(a \longmapsto n^{-\frac{1}{\kappa}}\left(T_{\psi(na)} - n \tau \psi(na) \right)\right)\xrightarrow[n\to\infty]{\sP_0^{(\alpha)}\text{-(d)}} c_1 \mathcal{S}^{\kappa}.
\]
This limit handles the first square bracket.

For the same reasons, the $\rho_{k+1} - \rho_k$ satisfy a central limit theorem, and, since $\kappa < 2$,
\[
\left(t \longmapsto n^{-\frac{1}{\kappa}}\left(\rho_{\lfloor nt \rfloor} - n \tau \langle v, \hat u \rangle t \right)\right)\xrightarrow[n\to\infty]{} 0.
\]
Thanks to Theorem 13.7.1 of \cite{whitt2002limits}: in $J_1$-topology,
\[
\left(a \longmapsto n^{-\frac{1}{\kappa}}\left(\psi(na) - n (\tau \langle v, \hat u \rangle)^{-1} a \right)\right)\xrightarrow[n\to\infty]{} 0.
\]
This limit handles the second square bracket.

All together, it implies:
$\left(a \longmapsto n^{-\frac{1}{\kappa}}\left(T_{\psi(na)} - n \langle v, \hat u \rangle^{-1} a \right)\right) \xrightarrow[n\to\infty]{\sP_0^{(\alpha)}\text{-(d)}} c_1 \mathcal{S}^{\kappa};$
and obviously, replacing $T_{\psi(na)}$ by $T_{\psi(na)-1}$ keeps this limit unchanged. As a consequence: in $J_1$-topology,
\[
\left(a \longmapsto n^{-\frac{1}{\kappa}}\left(\Delta_{na} - n \langle v, \hat u \rangle^{-1} a \right)\right)\xrightarrow[n\to\infty]{\sP_0^{(\alpha)}\text{-(d)}} c_2 \mathcal{S}^{\kappa}.
\]

\paragraph{Step 3: Limit theorem for $S$ and $\langle X, \hat u \rangle$.} 
Define, for $t \pg 0$:
$$S_t = \sup\ens{\langle X_k, \hat u \rangle, k \pp t}.$$

Note that $\bar{X}$ is the c\`adl\`ag inverse of $\Delta$. By Theorem 13.7.1 of \cite{whitt2002limits},
\[
\left(t \longmapsto n^{-\frac{1}{\kappa}}\left(S_{nt} - n \langle v, \hat u \rangle t \right)\right)\xrightarrow[n\to\infty]{\sP_0^{(\alpha)}\text{-(d)}} - c_3 \mathcal{S}^{\kappa}.
\]

But, thanks to condition $\mathbf{(T)}$, $\bar{X}$ and $\langle X, \hat u \rangle$ are very close to one another: if a final time $t^* > 0$ is fixed, since there are at most $nt^*$ renewals before time $nt^*$, for every $\epsilon > 0$, we have:
\begin{align*}
\sP_0\left(\sup_{t \pp nt^*} (S_t - \langle X_t, \hat u \rangle) > \epsilon n^{1/\kappa}\right)
& \pp \sum_{i=1}^{nt^*} \sP_0\left(\sup_{T_i \pp t \pp T_{i+1}} (S_t - \langle X_t, \hat u \rangle) > \epsilon n^{1/\kappa} \right) \\
& \pp \sum_{i=1}^{nt^*} \sP_0\left(\| X_{T_{i+1}} - X_{T_i} \| >  \epsilon n^{1/\kappa} \right) \\
& \pp C nt^* e^{-c \epsilon n^{1/\kappa}} \xrightarrow[n\to\infty]{} 0,
\end{align*}
so that, by Slutsky's lemma:
\[
\left(t \longmapsto n^{-\frac{1}{\kappa}}\left(\langle {X}_{nt}, \hat u \rangle - n \langle v, \hat u \rangle t \right)\right)\xrightarrow[n\to\infty]{\sP_0^{(\alpha)}\text{-(d)}} - c_3 \mathcal{S}^{\kappa}.
\]

\paragraph{Step 4: limit theorem for $X$.} We have got a one-dimensional result, about the fluctuations of the projection of $X$ along a fixed direction; we would like to get a $d$-dimensional result about $X$ itself. In order to that, we must control the fluctuations of the projection of $X$ over the orthogonal subspace of~$v$, the asymptotic direction.

Thanks to the previous point, we know that, in the asymptotic direction~$v$:
\[
\left(t \longmapsto n^{-\frac{1}{\kappa}}\left(\langle {X}_{nt}, v \rangle - n \| v \| t \right) v\right)\xrightarrow[n\to\infty]{\sP_0^{(\alpha)}\text{-(d)}} - c_3 \mathcal{S}^{\kappa} v.
\]

Let us now examine what happens in the orthogonal space of $v$. Let $\mathrm{pr}$ be the orthogonal projection along~$v$ over $v^\perp$, the orthogonal subspace of~$v$. For every $n \in \NN$, denote:
$k(n) = \sup\ens{k \pg 0 / T_k \pp n}$
the number of regenerations which occurred before time~$n$.

For every $n \in \NN$ and $t > 0$, we can decompose:
$$ \mathrm{pr}\left( \overline{X_{\lfloor nt \rfloor}} \right) = \mathrm{pr}\left(  \overline{X_{\lfloor nt \rfloor} - X_{T_{k(nt)}}} \right) + \mathrm{pr}\left( \overline{X_{T_{k(nt)}}} \right)$$
(in the sequel, overlining a random variable indicates that it is centered: $\overline{X} = X - \sE_0[X]$).

We have:
$$ \left( t \longrightarrow n^{-\frac{1}{\kappa}} \mathrm{pr}\left(\overline{ X_{\lfloor nt \rfloor} - X_{T_{k(nt)}} }\right) \right) \xrightarrow[n\to\infty]{\sP_0^{(\alpha)}} 0. $$
Indeed, for every $\epsilon > 0$ and $t^* > 0$, since there are at most $nt^*$ renewal times before time~$nt^*$, we have, by Markov's inequality:
\begin{align*}
&\sP_0 \left( \sup_{0 \pp t \pp t^*} \| \overline{X_{\lfloor nt \rfloor} - X_{T_{k(nt)}}}\| > \epsilon n^{\frac{1}{\kappa}} \right)
\pp \epsilon^{-2} n^{-\frac{2}{\kappa}} \sE_0 \left[ \sup_{0 \pp t \pp t^*} \| \overline{X_{\lfloor nt \rfloor} - X_{T_{k(nt)}}} \|^2 \right] \\*
\qquad &\pp 4 \epsilon^{-2} n^{-\frac{2}{\kappa}} \left( \sE_0 \left[ \sup_{0 \pp t \pp t^*} \| X_{\lfloor nt \rfloor} - X_{T_{k(nt)}} \|^2 \right] + \sE_0 \left[ \sup_{0 \pp t \pp t^*} \| X_{\lfloor nt \rfloor} - X_{T_{k(nt)}} \| \right]^2 \right) \\
\qquad &\pp 4 \epsilon^{-2} n^{-\frac{2}{\kappa}} \sum_{i=0}^{nt^*} \left( \sE_0 \left[ \sup_{T_i \pp nt \pp T_{i+1}} \| X_{\lfloor nt \rfloor} - X_{T_i} \|^2 \right] + \sE_0 \left[ \sup_{T_i \pp nt \pp T_{i+1}} \| X_{\lfloor nt \rfloor} - X_{T_i} \| \right]^2 \right) \\
\qquad &\pp 4 \epsilon^{-2} n^{-\frac{2}{\kappa}} \sum_{i=0}^{nt^*} \left( \sE_0 \left[ \sup_{T_i \pp k < T_{i+1}} \| X_k - X_{T_i} \|^2 \right] + \sE_0 \left[ \sup_{T_i \pp k < T_{i+1}} \| X_k - X_{T_i} \| \right]^2 \right) \\
\qquad &\pp C \epsilon^{-2} t^* n^{1-\frac{2}{\kappa}} \xrightarrow[n\to\infty]{} 0,
\end{align*}
by condition~$\mathbf{(T)}$.

Since the $X_{T_{k+1}} - X_{T_k}$ are i.i.d. and have exponential moments (because of condition~$\mathbf{(T)}$), we also have:
\[
\left(t \longmapsto n^{-\frac{1}{\kappa}}\mathrm{pr}(\overline{X_{T_{\lfloor nt \rfloor}}}) \right)\xrightarrow[n\to\infty]{\sP_0^{(\alpha)}\text{-(d)}} 0.
\]
Because $T_1$ is integrable, the $T_n$ satisfy a law of large number, and, as above, by inversion, it is also the case of the $k(n)$. By composition:
\[
\left(t \longmapsto n^{-\frac{1}{\kappa}}\mathrm{pr}(\overline{X_{T_{k(nt)}}}) \right)\xrightarrow[n\to\infty]{\sP_0^{(\alpha)}\text{-(d)}} 0.
\]

This concludes the proof.
\end{proof}

\section*{Appendix}

\begin{lemma}
\label{appendix}
Let $N$ be a random variable with geometric distribution with parameter~$p$. Then, for every $\beta > 0$, there exists a constant $C_\beta < \infty$ such that:
$$ \EE[N^\beta] < C_\beta p^{-\beta}.$$
\end{lemma}

\begin{proof}
It is well known that, for every non-negative integer $k \pg 0$:
$$ \EE[N (N-1) \ldots (N-k+1)] = \frac{k! (1-p)^{k-1}}{p^k}, $$
As the polynomial $X^k$ is a linear combination of the $X (X-1) \ldots (X-j+1)$, $0 \pp j \pp k$, this equality directly implies the existence of $C_k$.

We now show that, for every $k \in \ZZ_+$, the function:
$$ f_k: s \longmapsto p^{k+s} \EE[N^{k+s}], $$
defined for $s \in (0,1)$, is bounded. Let $k$ be a non-negative integer.

An easy application of the dominated convergence theorem shows that $f_k$ is twice differentiable, and that for every $s \in (0,1)$:
$$ f''_k(s) = \EE[(\ln(pN))^2 (pN)^{k+s}] \pg 0;$$
so $f_k$ is convex.

For every $s \in (0,1)$, by convexity:
\begin{align*}
f_k(s) &\pp (1-s) f_k(0) + s f_k(1) \\
&\pp (1-s) \EE[(pN)^{k}] + s \EE[(pN)^{k+1}] \\
&\pp (1-s) C_k + s C_{k+1} \\
&\pp C_k \vee C_{k+1}.
\end{align*}

It thus suffices to define  $C_\beta = C_{\lfloor \beta \rfloor} \vee C_{\lceil \beta \rceil}$ to prove the expected bound.
\end{proof}

\medskip

\textbf{Acknowledgements.} This work has been done during the preparation of my PhD under the supervision of Christophe Sabot. It is a pleasure for me to thanks him for valuable advice and constant support.

\bibliographystyle{plain}
\bibliography{subdiffusiveRWDE_bib}

\end{document}